\documentclass[10pt, twoside]{amsart}
\usepackage{pb-diagram}
\newcount\VOL\VOL=3\newcount\YEAR\YEAR=1998
\def\firstpage{177}\def\lastpage{205}
\def\received{March 22, 1998}\def\revised{September 21, 1998}
\def\communicated{Peter Schneider}
\font\eightrm=cmr8
\font\caps=cmcsc10                    % Theorem, Lemma etc
\font\Caps=cccsc10 scaled \magstep1   % Title
\font\scaps=cmcsc8
\makeatletter
\@specialpagefalse\headheight=8.5pt
\def\DocMath{{\def\th{\thinspace}\scaps Doc.\th Math.\th J.\th DMV}}
\renewcommand{\@oddfoot}{\hfill\scaps Documenta Mathematica 
    \number\VOL\  (\number\YEAR) \number\firstpage--\lastpage\hfill}
\renewcommand{\@evenfoot}{\ifnum\thepage>\lastpage\hfill\scaps
    Documenta Mathematica \number\VOL\  (\number\YEAR)\hfill\else\@oddfoot\fi}%

\renewcommand{\@evenhead}{%
    \ifnum\thepage>\lastpage\rlap{\thepage}\hfill%
    \else\rlap{\thepage}\slshape\leftmark\hfill\caps\SAuthor\hfill\fi}%
\renewcommand{\@oddhead}{%
    \ifnum\thepage=\firstpage{\DocMath\hfill\llap{\thepage}}%
    \else{\slshape\rightmark}\hfill\caps\STitle\hfill\llap{\thepage}\fi}%
\makeatother

\def\TSkip{\bigskip}
\newbox\TheTitle{\obeylines\gdef\GetTitle #1
\ShortTitle  #2
\SubTitle    #3
\Author      #4
\ShortAuthor #5
\EndTitle
{\setbox\TheTitle=\vbox{\baselineskip=20pt\let\par=\cr\obeylines%
\halign{\centerline{\Caps##}\cr\noalign{\medskip}\cr#1\cr}}%
        \copy\TheTitle\TSkip\TSkip%
\def\next{#2}\ifx\next\empty\gdef\STitle{#1}\else\gdef\STitle{#2}\fi%
\def\next{#3}\ifx\next\empty%
    \else\setbox\TheTitle=\vbox{\baselineskip=20pt\let\par=\cr\obeylines%
    \halign{\centerline{\caps##} #3\cr}}\copy\TheTitle\TSkip\TSkip\fi%
\centerline{\caps #4}\TSkip\TSkip%
\def\next{#5}\ifx\next\empty\gdef\SAuthor{#4}\else\gdef\SAuthor{#5}\fi%
\ifx\received\empty\relax
    \else\centerline{\eightrm Received: \received}\fi%
\ifx\revised\empty\TSkip%
    \else\centerline{\eightrm Revised: \revised}\TSkip\fi%
\ifx\communicated\empty\relax
    \else\centerline{\eightrm Communicated by \communicated}\fi\TSkip\TSkip%
\catcode'015=5}}\def\Title{\obeylines\GetTitle}
\def\Abstract{\begingroup\narrower
    \parskip=\medskipamount\parindent=0pt{\caps Abstract. }}
\def\EndAbstract{\par\endgroup\TSkip\TSkip}
\newbox\TheAdd\def\Addresses{\vfill\copy\TheAdd\vfill
    \ifodd\number\lastpage\vfill\eject\phantom{.}\vfill\eject\fi}
{\obeylines\gdef\GetAddress #1
\Address #2 
\Address #3
\Address #4
\EndAddress
{\def\xs{6truecm}
\setbox0=\vtop{{\obeylines\hsize=\xs#1}}\def\next{#2}
\ifx\next\empty % 1 address
     \setbox\TheAdd=\hbox to\hsize{\hfill\copy0\hfill}
\else\setbox1=\vtop{{\obeylines\hsize=\xs#2}}\def\next{#3}
\ifx\next\empty % 2 addresses
     \setbox\TheAdd=\hbox to\hsize{\hfill\copy0\hfill\copy1\hfill}
\else\setbox2=\vtop{{\obeylines\hsize=\xs#3}}\def\next{#4}
\ifx\next\empty\ % 3 addresses
     \setbox\TheAdd=\vtop{\hbox to\hsize{\hfill\copy0\hfill\copy1\hfill}
                \vskip20pt\hbox to\hsize{\hfill\copy2\hfill}}
\else\setbox3=\vtop{{\obeylines\hsize=\xs#4}}
     \setbox\TheAdd=\vtop{\hbox to\hsize{\hfill\copy0\hfill\copy1\hfill}
                \vskip20pt\hbox to\hsize{\hfill\copy2\hfill\copy3\hfill}}
\fi\fi\fi\catcode'015=5}}\gdef\Address{\obeylines\GetAddress}

\newtheorem*{theorem}{Theorem}
\newtheorem*{lemma}{Lemma}
\newtheorem*{keylemma}{Key Lemma}
\newtheorem*{proposition}{Proposition}
\newtheorem*{corollary}{Corollary}

\newcommand{\nobf}{ }
\theoremstyle{definition}

\theoremstyle{remark}
\newtheorem{remark}{Remark}

\numberwithin{equation}{subsection}
\numberwithin{remark}{section}

\newcommand{\N}{\mathbf{N}}
\newcommand{\Z}{\mathbf{Z}}

\newcommand{\ko}{\: , \;}

\newcommand{\myspace}{\makebox[1cm]{ }}
\newcommand{\mybigspace}{\makebox[2cm]{ }}
\newcommand{\ul}{\underline}

\newcommand{\bt}{\bullet}

\newcommand{\ra}{\rightarrow}
\newcommand{\la}{\leftarrow}

\newcommand{\arr}[1]{\stackrel{#1}{\rightarrow}}
\newcommand{\longarr}[1]{\xrightarrow{#1}}

\newcommand{\iso}{\stackrel{_\sim}{\rightarrow}}
\newcommand{\liso}{\stackrel{_\sim}{\leftarrow}}
\newcommand{\larr}[1]{\stackrel{#1}{\leftarrow}}

\renewcommand{\mod}{\operatorname{mod}\nolimits}
\newcommand{\Mod}{\operatorname{Mod}\nolimits}
\newcommand{\proj}{\operatorname{proj}\nolimits}
\newcommand{\coh}{\operatorname{coh}\nolimits}
\newcommand{\id}{\mathbf{1}}
\newcommand{\R}{\mathbf{R}}
\renewcommand{\L}{\mathbf{L}}

\newcommand{\ten}{\otimes}

\newcommand{\holim}{\operatorname{holim}}
\newcommand{\lid}{\varinjlim}
\newcommand{\lii}{\varprojlim}
\newcommand{\can}{\operatorname{can}}

\renewcommand{\H}[1]{{H}^{#1}}

\newcommand{\Hs}{{H}^*}

\newcommand{\HC}[1]{{HC}_{#1}\,}

\newcommand{\HCs}{{HC}_*}
\newcommand{\HNs}{{HC}^-_*}

\newcommand{\HCmixs}{{HC}_{mix,*}}
\newcommand{\HNmixs}{{HC}^-_{mix,*}}
\newcommand{\HPmixs}{{HC}^{per}_{mix,*}}
\newcommand{\HCders}{{HC}^{\mbox{\rm\scriptsize der}}_*}
\newcommand{\HCMcCs}{{HC}^{\mbox{\rm\scriptsize McC}}_*}
\newcommand{\HCMcC}[1]{{HC}^{\mbox{\rm\scriptsize McC}}_{#1}\,}

\newcommand{\ca}{{\mathcal A}}
\newcommand{\cb}{{\mathcal B}}
\newcommand{\cc}{{\mathcal C}}
\newcommand{\cd}{{\mathcal D}}

\newcommand{\ch}{{\mathcal H}}
\newcommand{\ci}{{\mathcal I}}

\newcommand{\cl}{{\mathcal L}}
\newcommand{\cm}{{\mathcal M}}

\newcommand{\co}{{\mathcal O}}

\newcommand{\cR}{{\mathcal R}}

\newcommand{\cs}{{\mathcal S}}
\newcommand{\ct}{{\mathcal T}}

\newcommand{\cx}{{\mathcal X}}
\newcommand{\eps}{\varepsilon}

\renewcommand{\phi}{\varphi}
\newcommand{\La}{\Lambda}

\newcommand{\Ga}{\Gamma}

\newcommand{\HOm}[3]{\operatorname{Hom}\nolimits_{#1}(#2,#3)}
\newcommand{\HOM}[3]{\ch om^\bullet_{\,#1}\,(#2,#3)}

\newcommand{\Cb}{\mathbf{C}^b}

\newcommand{\cone}{\operatorname{Cone}\nolimits}

\newcommand{\vecbundle}{\operatorname{vec}\nolimits}
\newcommand{\spec}{\operatorname{Spec}\nolimits}
\newcommand{\qcoh}{\operatorname{Qcoh}\nolimits}

\newcommand{\dmix}{\cd\cm ix\,}
\newcommand{\mormix}{\cm or \cm ix\,}
\newcommand{\dmormix}{\cd\mormix}

\newcommand{\per}{\operatorname{per}}
\newcommand{\on}{\mbox{ on }}
\newcommand{\strper}{\operatorname{strper}}
\newcommand{\flatper}{\operatorname{flatper}}
\newcommand{\iper}{\operatorname{iper}}
\newcommand{\tiper}{\widetilde{\operatorname{iper}}}
\newcommand{\roundper}{\wp er\,}

\newcommand{\Tot}{\operatorname{Tot}}
\newcommand{\hTot}{\widehat{\Tot}\,}
\begin{document} 
%%%% ------------- fill in your data below this line  -------------------

\Title
On the Cyclic Homology 
of Ringed Spaces and  Schemes
\ShortTitle
Cyclic homology
\SubTitle{ }   
\Author
Bernhard Keller
\ShortAuthor
\EndTitle

\Abstract   We prove that the cyclic homology of a scheme 
with an ample
line bundle coincides with the cyclic homology of its category
of algebraic vector bundles. As a byproduct of the proof, we
obtain a new construction of the Chern character of a perfect
complex on a ringed space.


\EndAbstract

\Address 
Bernhard Keller
UFR de Math\'ematiques
Universit\'e Paris 7
Institut Math\'ematique de Jussieu
UMR 9994 du CNRS
Case 7012
2, place Jussieu
75251 Paris Cedex 05
France
keller@math.jussieu.fr
www.math.jussieu.fr/$\tilde{\;\;\;}$keller

\Address
\Address
\Address
\EndAddress

%\begin{document}
%\input{montrere}

%\title{On the cylic homology of ringed spaces and
%schemes}

%\author{Bernhard Keller}
%\address{
%UFR de Math\'ematiques\\
%Universit\'e Paris 7\\
%Institut Math\'ematique de Jussieu\\
%UMR 9994 du CNRS\\
%Case 7012\\
%2, place Jussieu\\
%75251 Paris Cedex 05\\
%France
%}

%\email{
%\begin{minipage}[t]{5cm}
%keller@math.jussieu.fr \\
%www.math.jussieu.fr/$\tilde{\;\;\;}$keller
%\end{minipage}
%}

%\subjclass{16E40 (Primary), 18E30, 14F05 (Secondary)}
%\date{February 21, 1998}
%\keywords{Cyclic homology, Chern character, 
%Ringed space, Scheme, Perfect complex, Derived category}

%\begin{abstract}

%\end{abstract}

%\maketitle

%\tableofcontents

\section{Introduction}

\subsection{The main theorem.} Let $k$ be a field and $X$ a 
scheme over $k$ which admits an ample
line bundle (e.g. a quasi-projective variety). 
Let $\vecbundle(X)$ denote the category of algebraic 
vector bundles
on $X$. We view $\vecbundle(X)$ as an exact category in the sense
of Quillen \cite{Quillen73}: By definition, 
a short sequence of vector
bundles is admissible exact iff it is exact in the category of
sheaves on $X$. Moreover, the category $\vecbundle(X)$ 
is $k$-linear,
i.e. it is additive and its morphism sets are 
$k$-vector spaces such
that the composition is bilinear. In \cite{Keller97}, we 
have defined,
for each $k$-linear exact  category $\ca$, a cyclic homology
theory $\HCders(\ca)$. The superscript $der$ indicates that
the definition is modeled on that of the derived
category of $\ca$. In [loc. cit.] it was denoted by
$\HCs(\ca)$.
%Here and in the sequel,
%everything we say about `cyclic homology' also 
%applies mutatis mutandis
%to its variants: Hochschild homology as
%well as negative and periodic cyclic homology.
As announced in [loc. cit.],
in this article, we will show that the cyclic
homology of the scheme $X$ coincides with the cyclic homology
of the $k$-linear exact category $\vecbundle(X)$: 
There is a canonical isomorphism (cf. Corollary \ref{agreement})
\begin{equation}
\HCs(X) \iso \HCders(\vecbundle(X)). \label{gen_agreement}
\end{equation}
The definition of the cyclic homology of a scheme is
an important technical point which will be discussed below
in \ref{cychomschemes}. 
Note that by definition \cite[Par.~7]{Quillen73}, there is an
analogous isomorphism in $K$-theory.

\subsection{Motivation.}
Our motivation for proving 
the isomorphism \ref{gen_agreement} is
twofold: Firstly, it allows the computation of $\HCs(X)$ for
some non-trivial examples. 
Indeed, suppose that $k$ is algebraically
closed and that $X$ is a smooth projective algebraic variety.
Suppose moreover that $X$ admits a {\em tilting bundle}, i.e. 
a vector bundle without higher selfextensions whose direct
summands generate the bounded derived category of the category
of coherent sheaves on $X$. Examples of varieties satisfying
these hypotheses are projective spaces,
Grassmannians, and smooth quadrics 
\cite{Beilinson78}, 
\cite{Kapranov83}, \cite{Kapranov86}, \cite{Kapranov88inv}. 
In \ref{tiltingbundles},
we deduce from \ref{gen_agreement}
that for such a variety,
the Chern character induces an isomorphism
\[
K_0 X \ten_\Z \HCs k \iso \HCs X.
\]
Here the left hand side is explicitly
known since the group $K_0 X$ is free and admits a basis
consisting of the classes of the pairwise non-isomorphic
indecomposable direct summands of the tilting bundle.
Cyclic homology of projective spaces was first computed
by Beckmann \cite{Beckmann91} using a different method. 

Our second motivation for proving the isomorphism 
\ref{gen_agreement} is that it provides further justification
for the definition of $\HCders$. Indeed, there is
a `competing' (and previous) definition
of cyclic homology for $k$-linear exact categories
due to R.~McCarthy \cite{McCarthy94}. Let us denote by
$\HCMcCs(\ca)$ the graded $k$-module which he associates with
$\ca$. McCarthy proved in [loc. cit.] a number of good
properties for $\HCMcCs$. The most fundamental of these is
the existence of an agreement isomorphism
\[
\HCs (A) \iso \HCMcCs (\proj(A)) \ko
\]
where $A$ is a $k$-algebra and $\proj(A)$, the category of finitely
generated projective $A$-modules endowed with the
split exact sequences. In particular, if we take $A$ to be
commutative, we obtain the isomorphism
\[
\HCs(X) \iso \HCMcCs(\vecbundle(X))
\]
for all {\em affine schemes $X=\spec(A)$} 
(to identify the left hand side,
we use Weibel's isomorphism \cite{Weibel96} between the
cyclic homology of an affine scheme and the cyclic homology
of its coordinate algebra). Whereas for $\HCders$, 
this ismorphism
extends to more general schemes, this cannot be the case for
$\HCMcCs$. Indeed, for $n\geq 0$, the group $\H{n}(X, \co_X)$
occurs as a direct factor of $\HC{-n}(X)$. However, the
group $\HCMcC{-n}$ vanishes for $n>0$ by its very definition.

\subsection{Generalization, Chern character.}
Our proof of the isomorphism \ref{gen_agreement} actually yields
a more general statement: Let $X$ be a quasi-compact separated
scheme over $k$. Denote by $\per X$ the pair formed by the 
category of perfect sheaves (\ref{presheafperfect}) on $X$
and its full subcategory of acyclic perfect sheaves. The
pair $\per X$ is a localization pair in the sense of
\cite[2.4]{Keller97} and its cyclic homology 
$\HCs(\per X)$ has been defined in [loc. cit.]. We will show
(\ref{agreement}) that there is a canonical isomorphism
\begin{equation}
\HCs(X) \iso \HCs(\per X).\label{perf_agreement}
\end{equation}
If $X$ admits an ample line bundle, we have
an isomorphism
\[
\HCders(\vecbundle(X)) \iso \HCs(\per X)
\]
so that the isomorphism \ref{gen_agreement} results as
a special case. 

The first step in the proof of \ref{perf_agreement} will be
to construct a map
\[
\HCs(\per X) \ra \HCs(X).
\]
This construction will be carried out in \ref{charclasses}
for an arbitrary 
topological space $X$ endowed with a sheaf of (possibly
non-commutative) $k$-algebras. 
As a byproduct, we therefore
obtain a new construction of the Chern character of a
perfect complex $P$. Indeed, the complex $P$ yields a
functor between localization pairs
\[
?\ten_k P : \per \mbox{pt} \ra \per X
\]
and hence a map
\[
\HCs(\per \mbox{pt}) \ra \HCs(\per X) \ra \HCs(X).
\]
The image of the class
\[
ch([k]) \in \HCs(\per \mbox{pt}) = \HCs(k)
\]
under this map is the value of the Chern character at the
class of $P$. An analogous construction works for the
other variants of cyclic homology, in particular for
negative cyclic homology. The first construction of a
Chern character for perfect complexes is
due to Bressler--Nest--Tsygan, who needed it in
their proof \cite{BresslerNestTsygan97} of 
Schapira-Schneiders' conjecture \cite{SchapiraSchneiders94}.
They even construct a generalized Chern character defined
on all higher $K$-groups. Several other constructions
of a classical Chern character are due to B.~Tsygan
(unpublished).

\subsection{Cyclic homology of schemes.}
\label{cychomschemes}
Let $k$ be a commutative ring and $X$ a scheme over $k$. 
The cyclic homology of $X$ was first defined by 
Loday \cite{Loday86}: He sheafified the classical bicomplex
to obtain a complex of sheaves $CC(\co_X)$. He then defined
the cyclic homology of $X$ to be the hypercohomology of the
(total complex of) $CC(\co_X)$. Similarly for the different
variants of cyclic homology. There arise three problems:
\begin{itemize}
\item[(1)] The complex $CC(\co_X)$ is unbounded to the left.
So there are (at least) two non-equivalent possibilities to
define its hypercohomology: should one take Cartan-Eilenberg
hypercohomology (cf. \cite{Weibel96}) or derived functor
cohomology in the sense of Spaltenstein \cite{Spaltenstein88} ?
\item[(2)] Is the cyclic homology of an affine scheme isomorphic
to the cyclic homology of its coordinate ring ?
\item[(3)] If a morphism of schemes induces an isomorphism in
Hochschild homology, does it always induce an isomorphism in
cyclic homology ?
\end{itemize}
Problem (1) is related to the fact that in a category of
sheaves, products are not exact in general. We refer to
\cite{Weibel96} for a discussion of this issue.

In the case of a noetherian scheme of finite dimension,
Beckmann \cite{Beckmann91} and Weibel-Geller \cite{WeibelGeller91}
gave a positive answer to (2) using Cartan-Eilenberg 
hypercohomology. By proving the existence of an SBI-sequence
linking cyclic homology and Hochschild homology they also
settled (3) for this class of schemes, whose Hochschild
homology vanishes in all sufficiently negative degrees.
Again using Cartan-Eilenberg hypercohomology, Weibel gave
a positive answer to (2) in the general case in \cite{Weibel96}.
There, he also showed that cyclic homology is a homology
theory on the category of quasi-compact quasi-separated
schemes. Problem (3) remained open. 

We will show in \ref{sheaveswith} that Cartan-Eilenberg
hypercohomology agrees with Spaltenstein's derived functor
hypercohomology on all complexes with quasi-coherent homology
if $X$ is quasi-compact and separated. Since $CC(\co_X)$ has
quasi-coherent homology
\cite{WeibelGeller91}, this shows that problem (1) does
not matter for such schemes. As a byproduct of \ref{sheaveswith},
we deduce in \ref{quasicoherentsheaves} a
(partially) new proof of Boekstedt-Neeman's theorem
\cite{BoekstedtNeeman93} which states that for a quasi-compact
separated scheme $X$, the unbounded derived category of
quasi-coherent sheaves on $X$ is equivalent to the full subcategory
of the unbounded derived category of all $\co_X$-modules whose
objects are the complexes with quasi-coherent homology.  A different
proof of this was given by Alonso--Jerem\'{\i}as--Lipman in
\cite[Prop. 1.3]{AlonsoJeremiasLipman97}.

In order to get rid of problem (3), we will slightly
modify Loday's definition: Using sheaves of mixed complexes
as introduced by Weibel \cite{Weibel97} we will show that
the image of the Hochschild complex $C(\co_X)$ under
the derived global section functor is canonically a
mixed complex $M(X)$. The {\em mixed cyclic homology} of
$X$ will then be defined as the cyclic homology of $M(X)$. 
For the mixed cyclic homology groups, the answer 
to (2) is positive thanks to the corresponding theorem
in Hochschild homology due to Weibel--Geller \cite{WeibelGeller91};
the answer to (3) is positive thanks to the
definition. The mixed cyclic homology groups coincide with Loday's
groups if the derived global section functor commutes
with infinite sums. This is the case for quasi-compact
separated schemes as we show in \ref{proofofcorollary}.

\subsection{Organization of the article.} 
In section 2, we recall the mixed complex of an algebra
and define the mixed complex $M(X,\ca)$ of a ringed
space $(X,\ca)$. In section 3, we recall the definition
of the mixed complex associated with a localization
pair and give a `sheafifiable' 
description of the Chern character of
a perfect complex over an algebra. 
In section 4, we construct a morphism from the
mixed complex associated with the category of perfect
complexes on $(X,\ca)$ to the mixed complex $M(X,\ca)$.
We use it to construct the Chern character of a perfect
complex on $(X,\ca)$. In section 5, we state and prove
the main theorem and apply it to the computation of
the cyclic homology of smooth projective varieties
admitting a tilting bundle.
In appendix A, we prove that Cartan-Eilenberg
hypercohomology coincides with derived functor cohomology
for (unbounded) complexes with quasi-coherent homology 
on quasi-compact separated schemes. In appendix B, we apply
this to give a (partially) new proof of a theorem
of Boekstedt-Neeman \cite{BoekstedtNeeman93}.

\subsection{Acknowledgment.} The author thanks the referee
for his suggestions, which helped to 
make this article more readable.
\section{Homology theories for ringed spaces}
Let $k$ be a field, 
$X$ a topological space, and $\ca$ a sheaf of 
$k$-algebras on $X$. In this section, we consider
the possible definitions of the cyclic homology of
$(X,\ca)$. In \ref{Hochschild} we recall the definition suggested
by Loday \cite{Loday86}. In \ref{Mixed}, we point out
that with this definition, it is not clear that a 
morphism inducing isomorphisms in Hochschild homology
also does so in cyclic homology and its variants.
This is our main reason for introducing the
`mixed homologies'. These also have the advantage
of allowing a unified and simultaneous
treatment of all the different homology theories.
For the sequel, the two fundamental invariants are
the mixed complex of sheaves $M(\ca)$ and its
image $M(X,\ca)=\R \Gamma(X, M(\ca))$ under the derived global
section functor. Both are canonical up to quasi-isomorphism
and are therefore viewed as objects of the corresponding
mixed derived categories. In the
case of a point and a sheaf given by an algebra $A$, 
these complexes specialize to the mixed complex
$M(A)$ associated with the algebra. The mixed
cyclic homology $\HCmixs(\ca)$ is defined to be
the cyclic homology of the mixed complex $M(X,\ca)$.

\subsection{Hochschild and cyclic homologies.}
\label{Hochschild}
Following a suggestion by Loday \cite{Loday86},
the {\em Hochschild complex} $C(\ca)$,
and the bicomplexes $CC(\ca)$, $CC^-(\ca)$, and 
$CC^{per}(\ca)$ are 
defined in \cite[4.1]{BresslerNestTsygan97} by
composing the classical constructions
(cf. \cite{Loday92}, for example) with
sheafification. The {\em Hochschild homology}, 
{\em cyclic homology}  \ldots of $\ca$ are then obtained as the
homologies of the complexes 
\[
\R\Gamma(X, C(\ca)) \ko \R\Gamma(X, CC(\ca))\ko \ldots 
\]
where $\R\Gamma(X,?)$ is the total right derived functor
in the sense of Spaltenstein \cite{Spaltenstein88} of the
global section functor.

\subsection{Mixed cyclic homologies.}
\label{Mixed}
Suppose that $f:(X,\ca)\ra (Y,\cb)$ is a morphism
of spaces with sheaves of $k$-algebras 
inducing isomorphisms in Hochschild homology.
With the above definitions, it does not seem to
follow that $f$ also induces isomorphisms
in cyclic homology, negative cyclic homology, and
periodic cyclic homology.
This is one of the reasons why we need to replace
the above definitions by slightly different variants 
defined in terms of the mixed complex associated
with $\ca$. This complex was introduced by
C.~Weibel in \cite{Weibel97}. However, the
`mixed homologies' we consider do not always
coincide with the ones of \cite{Weibel97} (cf. the end of this section).

Let us first recall the case of ordinary algebras:
For an algebra $A$, we denote by $M(A)$ the mapping
cone over the differential $1-t$ linking the first two
columns of the bicomplex $CC(A)$. We endow $M(A)$ with
the operator $B: M(A) \ra M(A)[-1]$ induced by
the norm map $N$ from the first to the second column
of the bicomplex. Then endowed with its differential $d$
and with the operator $B$ the complex $M(A)$
becomes a mixed complex in the sense of Kassel
\cite{Kassel87}, i.e. we have
\[
d^2=0 \ko B^2=0 \ko dB+Bd=0.
\]
The mixed complex $M(A)$ completely
determines the homology theories of $A$. Indeed,
we have a canonical quasi-isomorphism
\[
C(A) \ra M(A),
\]
which shows that Hochschild homology is determined by $M(A)$.
We also have canonical quasi-isomorphisms
\[
CC(A) \iso M(A) \ten_\Lambda^\L k \quad \ko \quad
CC^-(A) \iso \R \HOm{\Lambda}{k}{M(A)}
\]
where the right hand sides are defined by viewing
mixed complexes as objects of the mixed derived category,
i.e. differential graded (=dg) modules
over the dg algebra $\Lambda$
generated by an indeterminate $\eps$ of
chain degree $1$ with $\eps^2=0$ and $d\eps=0$
(cf. \cite{Kassel87}, \cite{Keller94}).
Finally, we have a quasi-isomorphism
\[
CC^{per}(A) \ra (\R\lii) P_k[-2n]\ten_\Lambda M(A)
\]
where $P_k$ is a cofibrant resolution 
(= `closed' resolution in the sense of \cite[7.4]{Keller98}
=`semi-free' resolution in the sense of 
\cite{AvramovFoxbyHalperinxy})
of the dg $\Lambda$-module $k$ and the transition map
$P_k[-2(n+1)] \ra P_k[-2n]$
comes from a chosen morphism of mixed complexes $P_k \ra P_k[2]$
which represents the canonical morphism $k \ra k[2]$
in the mixed derived category.
For example, one can take 
\[
P_k= \bigoplus_{i\in \N} \Lambda[2i]
\]
as a $\Lambda$-module endowed with the 
differential mapping the generator $1_{i}$
of $\Lambda[2i]$ to $\eps 1_{i-1}$.
The periodicity morphism then takes
$1_i$ to $1_{i-1}$ and $1_0$ to $0$.
Note that the functor
$\lii P_k[-2n]\ten_\Lambda ?$
is actually exact so that $\R\lii$
may be replaced by $\lii$ in the
above formula. 

Following Weibel \cite[Section 2]{Weibel97}
we sheafify this construction to obtain a mixed
complex of sheaves $M(\ca)$. We view it
as an object of the {\em mixed derived category 
$\dmix(X)$ of sheaves on $X$}, i.e. the derived category of
dg sheaves over the constant sheaf of
dg algebras with value $\Lambda$. The global
section functor induces a functor from mixed
complexes of sheaves to mixed complexes
of $k$-modules. By abuse of notation, the
total right derived functor of the induced
functor will still be denoted by $\R\Ga(X,?)$.
The {\em mixed complex of the ringed space $(X,\ca)$}
is defined as
\[
M(X,\ca)= \R\Ga(X,M(\ca)).
\]
The fact that the functor $\R\Ga(X,?)$ (and the mixed
derived category of sheaves) is well defined
is proved by adapting Spaltenstein's argument
of section 4 of \cite{Spaltenstein88}. 
Since the underlying complex of $k$-modules of
$M(\ca)$ is quasi-isomorphic to $C(\ca)$,
we have a canonical isomorphism
\[
HH_*(\ca) \iso H_* \R\Ga(X, M(\ca)).
\]
We define the `mixed variants'
\[
\HCmixs(\ca)\ko \HNmixs(\ca) \ko \HPmixs(\ca)
\]
of the homologies associated with $\ca$ by
applying the functors
\[
?\ten_\La^L k\ko 
\R\HOm{\Lambda}{k}{?} \quad\mbox{resp.} \quad
\R\lii P_k[-2n]\ten_\Lambda ?
\]
to $M(X,\ca)$ and taking homology.

These homology theories are slightly different from those of
Bressler--Nest--Tsygan \cite{BresslerNestTsygan97}, Weibel
\cite{Weibel96}, \cite{Weibel97}, and Beckmann \cite{Beckmann91}. 
We prove in \ref{proofofcorollary}
that mixed cyclic homology coincides with the cyclic homology defined
by Weibel if the global section functor $\R\Ga(X,?)$ commutes with
countable coproducts and that this is the case if $(X,\ca)$ is
a quasi-compact separated scheme.

For a closed subset $Z\subset X$,
we obtain versions with support in $Z$
by applying the corresponding functors to $\R\Ga_Z(X,M(\ca))$.

Now suppose that a morphism $(X,\ca)\ra (Y,\cb)$ induces
an isomorphism in $HH_*$. Then by definition, it induces
an isomorphism in the mixed derived category
\[
\R\Ga(X,M(\ca)) \la \R\Ga(Y,M(\cb))
\]
and thus in $\HCmixs$, $\HNmixs$, and $\HPmixs$.

%%%%%%%%%%%%%%%%%%%%%%%%%%%%%%%%%%%%%%%%%%%%%
%%%%%%%%%%%%%%%%%%%%%%%%%%%%%%%%%%%%%%%%%%%%%%%%%%%%%%%%%%%%%%%%%%%%%%%

\section{Homology theories for categories}

In this section, we recall the definition of the
cyclic homology (or rather: the mixed complex) 
of a localization pair from \cite{Keller97}. 
We apply this to give a description of the
Chern character of a perfect complex 
over an algebra $A$ (=sheaf of algebras over
a point). This description will later be
generalized to sheaves of algebras over
a general topological space.

A localization pair is a pair consisting of
a (small) differential graded $k$-category
and a full subcategory satisfying certain
additional assumptions. To define its mixed complex,
we proceed in three steps: In \ref{kcategories},
the classical
definition for algebras is generalized to small
$k$-categories following an idea of Mitchell's
\cite{Mitchell72}; then, in \ref{dgcategories},
we enrich our small $k$-categories
over the category of differential complexes, i.e. we
define the mixed complex of a differential graded
small $k$-category; by making this definition relative
we arrive, in \ref{pairs},
at the definition of the mixed complex of a
localization pair. For simplicity, we work only with
the Hochschild complex at first. 

We illustrate each
of the three stages by considering the respective
categories associated with a $k$-algebra $A$ :
the $k$-category $\proj(A)$ of finitely generated
projective $A$-modules, the differential graded
$k$-category $C^b(\proj(A))$ of bounded complexes
over $\proj(A)$, and finally the localization pair
formed by the category of all perfect complexes
over $A$ together with its full subcategory of
all acyclic perfect complexes. The three respective
mixed complexes are canonically quasi-isomorphic.
Thanks to this fact the mixed complex of an algebra
is seen to be functorial with respect to exact
functors between categories of perfect complexes.
This is the basis for our description of the
Chern character in \ref{charclasses}.

\subsection{$k$-categories.}
\label{kcategories}
Let $\cc$ be a small $k$-category, i.e. a small category
whose morphism spaces carry structures of $k$-modules
such that the composition maps are bilinear.
Following Mitchell \cite{Mitchell72} one defines
the Hochschild complex $C(\cc)$ to be the complex
whose $n$th component is
\begin{equation}
\label{directsum}
\coprod \cc(X_n, X_0) \ten \cc(X_{n-1}, X_n) \ten 
\cc(X_{n-2}, X_{n-1}) \ten \ldots \ten \cc(X_0, X_1) 
\end{equation}
where the sum runs over all sequences $X_0, \ldots, X_n$ of
objects of $\cc$. The differential is given by the
alternating sum of the face maps
\[
d_i (f_{n}, \ldots, f_i, f_{i-1}, \ldots, f_0)=
\left\{
\begin{array} {ll}
(f_{n}, \ldots, f_i f_{i-1}, \ldots f_0) & \mbox{ if $i>0$} \\
(-1)^{n} (f_0 f_{n}, \ldots, f_1) & \mbox{ if $i=0$}
\end{array}
\right.
\]
For example, suppose that $A$ is a $k$-algebra.
If we view $A$ as a category $\cc$ with one object, the
Hochschild complex $C(\cc)$ coincides with $C(A)$.
We have a canonical functor
\[
A \ra \proj A\ko
\]
where $\proj A$ denotes the category of finitely
generated projective $A$-modules. By a theorem
of McCarthy \cite[2.4.3]{McCarthy94}, this functor induces
a quasi-isomorphism 
\[
C(A) \ra C(\proj A).
\]

\subsection{Differential graded categories.}
\label{dgcategories}
Now suppose that the category $\cc$ is a differential
graded $k$-category. This means that $\cc$ is enriched
over the category of differential $\Z$-graded 
$k$-modules (=dg $k$-modules), i.e. each
space $\cc\,(X,Y)$ is a dg $k$-module and the
composition maps
\[
\cc\,(Y,Z)\ten_k \cc\,(X,Y) \ra \cc\,(X,Z)
\]
are morphisms of dg $k$-modules.
Then we obtain a  double complex whose columns are
the direct sums of (\ref{directsum}) and whose horizontal
differential is the alternating sum of the
face maps
\[
d_i (f_{n}, \ldots, f_i, f_{i-1}, \ldots, f_0)=
\left\{
\begin{array} {ll}
(f_{n}, \ldots, f_i f_{i-1}, \ldots f_0) & \mbox{ if $i>0$} \\
(-1)^{(n+\sigma)} (f_0 f_{n}, \ldots, f_1) & \mbox{ if $i=0$}
\end{array}
\right.
\]
where $\sigma= (\deg f_0) (\deg f_1 + \cdots + \deg f_{n-1})$.
The Hochschild complex $C(\cc)$ of the dg category $\cc$ is by
definition the (sum) total complex of this double
complex. 
The dg categories we will encounter are all obtained
as subcategories of a category $\mathbf{C}(\cx)$
of differential complexes over a $k$-linear
category $\cx$ (a $k$-linear category is a 
$k$-category which admits all finite direct sums). 
In this case, the dg structure
is given by the complex $\HOM{\cx}{X}{Y}$ associated
with two differential complexes $X$ and $Y$.

Hence if $A$ is a $k$-algebra,
the category $\Cb(\proj A)$ of bounded complexes
of finitely generated projective $A$-modules is
a dg category and the functor
\[
\proj A \ra \Cb(\proj A)
\]
mapping a module $P$ to the complex concentrated
in degree $0$ whose zero component is $P$ becomes
a dg functor if we consider $\proj A$ as a
dg category whose morphism spaces are concentrated
in degree 0. 
By \cite[lemma 1.2]{Keller98},
the functor $\proj A \ra \Cb(\proj A)$
induces a quasi-isomorphism
\[
C(\proj A) \ra C(\Cb(\proj A)).
\]

\subsection{Pairs of dg categories.}
\label{pairs}
Now suppose that
$\cc_0 \subset \cc_1$ are full subcategories
of a category of complexes $\mathbf{C}(\cx)$
over a small $k$-linear category $\cx$.
We define the Hochschild complex $C(\cc)$
of the pair $\cc: \cc_0 \subset \cc_1$
to be the cone over the morphism
\[
C(\cc_0) \ra C(\cc_1)
\]
induced by the inclusion (here both $\cc_0$
and $\cc_1$ are viewed as dg categories).
For example, let $A$ be a $k$-algebra. Recall
that a perfect complex over $A$ is a complex
of $A$-modules which is quasi-isomorphic to
a bounded complex of finitely generated projective
$A$-modules. Let $\per A$ denote the pair
of subcategories of the category of complexes
of $A$-modules formed by the category
$\per_1 A$ of perfect $A$-modules and its
full subcategory $\per_0 A$ of acyclic perfect
$A$-modules. Clearly we have a functor
$\proj A \ra \per A$, i.e. a commutative
diagram of dg categories
\[
\begin{diagram}
\node{0} \arrow{e} \arrow{s} 
\node{\per_0 A} \arrow{s} \\
\node{\proj A} \arrow{e}
\node{\per_1 A}
\end{diagram}
\]
This functor induces a quasi-isomorphism
\[
C(\proj A) \ra C(\per A)
\]
by theorem 2.4 b) of \cite{Keller97}.

\subsection{Mixed complexes and characteristic classes.}
\label{mixedcomcharclass}
In the preceding paragraph, we have worked with
the Hochschild complex, but it is easy to check that
everything we said carries over to the mixed complex
(\ref{Mixed}). The conclusion is then that
if $A$ is a $k$-algebra, we have
the following isomorphisms in the mixed derived category
\[
M(A) \iso M(\proj A) \iso M(\per A).
\]
This shows that $M(A)$ is functorial with respect to
morphisms of pairs $\per A \ra \per B$, i.e. functors
from perfect complexes over $A$ to perfect complexes
over $B$ which respect the dg structure and preserve
acyclicity. For example, if $P$ is a perfect complex
over $A$, we have the functor
\[
?\ten_k P : \per k \ra \per A
\]
which induces a morphism
\[
M(?\ten_k P) : M(\per k) \ra M(\per A)
\]
and hence a morphism
\[
M(P): M(k) \ra M(A).
\]
If we apply the functors $H_0$ resp. $\Hs\R\HOm{\Lambda}{k}{?}$
to this morphism we obtain morphisms
\[
HH_0(k) \ra HH_0(A) \quad\mbox{and}\quad
\HNmixs(k) \ra \HNmixs(A)
\]
which map the canonical classes in $HH_0(k)$
resp. $\HNmixs(k)=\HNs(k)$ to the Euler class resp. the
Chern character of the perfect complex $P$.

\section{Characteristic Classes for Ringed spaces}

Let $k$ be a field, 
$X$ a topological space, and $\ca$ a sheaf of 
$k$-algebras on $X$. In this section, we consider, 
for each open subset
$U$ of $X$, the localization pair of perfect complexes
on $U$ denoted by $\per \ca|_U$. The mixed complexes
$M(\per \ca|_U$) associated with these localization pairs 
are assembled into a sheaf of mixed complexes $M(\roundper \ca)$.
In \ref{presheafperfect}, we show that this sheaf is
quasi-isomorphic to the sheaf $M(\ca)$ of mixed complexes
associated with $\ca$. In \ref{charclasses}, this isomorphism
is used to construct the trace morphism
\[
\tau : M(\per \ca) \ra \R\Ga(X, M(\ca)).
\]
The construction of the characteristic classes of
a perfect complex is then achieved using the functoriality
of the mixed complex $M(\per \ca)$ with respect to
exact functors between localization pairs.

The main theorem (\ref{agreement}) will state that 
$\tau$ is invertible if $(X,\ca)$ is a quasi-compact
separated scheme.

\subsection{The presheaf of categories of perfect complexes.}
\label{presheafperfect}
Recall that a {\em strictly perfect complex} is a 
complex $P$ of $\ca$-modules such that each
point $x\in X$ admits an open neighbourhood $U$ such that
$P|_U$ is isomorphic to a bounded complex of
direct summands of finitely generated free 
$\ca|_U$-modules (note that such
modules have no reason to be projective objects in the 
category of $\ca|_U$-modules). A {\em perfect complex} is
a complex $P$ of $\ca$-modules such that each
point $x\in X$ admits an open neighbourhood
$U$ such that $P|_U$ is quasi-isomorphic to
a strictly perfect complex.

We denote by $\per \ca$ the pair formed by
the category of perfect complexes and its
full subcategory of acyclic perfect complexes.
For each open $U\subset X$, we denote by
$\per \ca|_U$ the corresponding pair of categories
of perfect $\ca|_U$-modules. Via the restriction
functors, the assignment $U \mapsto M(\per (\ca|_U))$
becomes a presheaf of mixed complexes on $X$.
We denote by $M(\roundper \ca)$ the corresponding
sheaf of mixed complexes.

For each open $U\subset X$, we have a canonical
functor
\[
\proj \ca(U) \ra \per \ca|_U \ko
\]
whence morphisms
\[
M(\ca(U)) \ra M(\proj \ca(U)) \ra M(\per\ca|_U)
\]
and a morphism of sheaves
\[
M(\ca) \ra M(\roundper \ca).
\]
\begin{keylemma} \label{keylemma} The above morphism
is a quasi-isomorphism
\end{keylemma}
\begin{remark} This is the analog in cyclic homology
of lemma 4.7.1 of \cite{BresslerNestTsygan97} (with the
same proof, as P.~Bressler has kindly informed me).
\end{remark}
\begin{proof} We will show that the morphism induces
quasi-isomorphisms in the stalks. Let $x\in X$. Clearly we have
an isomorphism
\[
M(\roundper \ca)_x \iso M(\lid \per \ca|_U)\ko
\]
where $U$ runs through the system of open neighbourhoods
of $x$. We will show that the canonical functor
\[
\lid \per \ca|_U \ra \per \ca_x
\]
induces a quasi-isomorphism in the mixed complexes.
For this, it is enough to show that it induces equivalences
in the associated triangulated categories, by \cite[2.4 b)]{Keller97}.
Now we have a commutative square
\[
\begin{diagram}
\node{\lid\per \ca|_U} \arrow{e} 
\node{\per \ca_x} 
\\
\node{\lid \strper \ca|_U} \arrow{e} \arrow{n}
\node{\strper \ca_x} \arrow{n}
\end{diagram}
\]
Here, we denote by $\strper$ the pair formed by the
category of strictly perfect complexes and its subcategory
of acyclic complexes. For an algebra $A$, we have
$\strper A = \Cb(\proj A)$ by definition.
It is easy to see that the two vertical arrows
induce equivalences in the triangulated categories,
and the bottom arrow is actually itself an equivalence
of categories. Indeed, we have the commutative square
\[
\begin{diagram}
\node{\lid \strper \ca|_U} \arrow{e} 
\node{\strper \ca_x}
\\
\node{\lid \strper \ca(U)} \arrow{n} \arrow{e}
\node{\strper\ca_x} \arrow{n}
\end{diagram}
\]
Here the right vertical arrow is the identity and
the left vertical arrow and the bottom arrow
are clearly equivalences.

The claim follows since the composition of the
morphism 
\[
M(\ca_x) \ra M(\lid \per \ca|_U)
\]
with the quasi-isomorphism $M(\lid \per \ca|_U) \ra M(\per \ca_x)$
is the canonical quasi-isomorphism $M(\ca_x) \ra M(\per\ca_x)$.

\end{proof}

\subsection{Characteristic classes.}
\label{charclasses}
By definition of $M(\roundper \ca)$ we have a morphism
of mixed complexes $M(\per \ca) \ra \Gamma(X, M(\roundper\ca))$.
By the key lemma (\ref{presheafperfect}), the canonical morphism
$M(\ca) \ra M(\roundper\ca)$ is invertible in the
mixed derived category. Thus we
can define the trace morphism 
\[
\tau : M(\per \ca) \ra \R\Gamma(X, M(\ca))
\]
by the following commutative diagram
\[
\begin{diagram}
\node{M(\per \ca)} \arrow{e} \arrow{s,l}{\tau} 
\node{\Gamma(X, M(\roundper \ca))} \arrow{s}
\\
\node{\R\Gamma(X, M(\ca))} \arrow{e}
\node{\R\Gamma(X, M(\roundper\ca))}
\end{diagram}
\]
Now let $P$ be a perfect complex.
It yields a functor
\[
?\ten_k P : \per k \ra \per \ca
\]
and hence a morphism in the mixed derived category
\[
M(k) \iso M(\per k) \xrightarrow{M(P)} 
M(\per \ca) \arr{\tau} \R\Gamma(X, M(\ca))=M(X,\ca).
\]
If we apply the functor $H_0$ resp. $\R\HOm{\Lambda}{k}{?}$
to this morphism, we obtain morphisms
\[
HH_0(k) \ra HH_0(\ca) \quad\mbox{resp.}\quad
\HNs(k)=\HNmixs(k) \ra \HNmixs(\ca)
\]
mapping the canonical classes to the Euler class respectively
to the Chern character of the perfect complex $P$.

\begin{remark} The trace morphism 
$\tau: M(\per \ca) \ra M(X,\ca)$ is
a quasi-isomor\-phism if $X$ is a point (by \ref{pairs}) or if
$(X,\ca)$ is a quasi-compact separated scheme
(by \ref{agreement} below).
\end{remark}

\begin{remark} (B.~Tsygan) Let $P$ be a perfect complex and 
$A=\HOM{X}{P}{P}$ the dg algebra of endomorphisms of $P$.
So if $P$ is fibrant (cf. \ref{terminology}),
then the $i$th homology of $A$ identifies with 
$\HOm{\cd X}{P}{P[i]}$. The dg category with
one object whose endomorphism algebra is $A$ naturally
embeds into $\per_1 \ca$ and we thus obtain a morphism
\[
M(A) \ra M(\per_1 \ca ) \ra M(\per \ca) \arr{\tau} \R\Ga(X,M(\ca))
\]
whose composition with the canonical map 
$M(k) \ra M(A)$ coincides with the morphism constructed above.
\end{remark}

\subsection{Variant with supports.}
\label{withsupports}
Let $Z\subset X$ be a closed subset. Let
$\per(\ca \on X)$ be the pair formed by the category
of perfect complexes acyclic off $Z$ and its full
subcategory of acyclic complexes. For each open
$U\subset X$ denote by $\per(\ca|_U \on Z)$ the
corresponding pair of categories of perfect
$\ca|_U$-modules. Via the restriction functors,
the assignment $U \mapsto M(\per(\ca|_U \on Z))$ 
becomes a presheaf of mixed complexes on $X$. We denote
by $M(\roundper(\ca\on Z))$ the corresponding sheaf
of mixed complexes. We claim that $M(\roundper(\ca\on Z))_x$
is acyclic for $x\not\in Z$. Indeed, if $U\subset X\setminus Z$
is an open neighbourhood of $x$, then by definition, the
inclusion
\[
\per_0(\ca|_U\on Z) \ra \per_1(\ca|_U\on Z)
\]
is the identity so that $M(\per(\ca|_U\on Z))$ is
nullhomotopic. It follows that the canonical
morphism $M(\roundper(\ca\on Z)) \ra M(\roundper\ca)$
uniquely factors through 
\[\R\Ga_Z M(\roundper\ca) \ra M(\roundper\ca)
\]
in $\dmix(X)$. Using the quasi-isomorphism 
$M(\ca) \ra M(\roundper\ca)$ we thus obtain a canonical morphism
$M(\roundper (\ca\on Z)) \ra \R\Ga_Z M(\ca)$ making
the following diagram commutative
\[
\begin{diagram}
\node{M(\roundper(\ca \on Z))} \arrow{e} \arrow{se} 
\node{\R\Ga_Z M (\roundper \ca)} \arrow{e} 
\node{M(\roundper\ca)} \\
\node{ }
\node{\R\Ga_Z M(\ca)} \arrow{n} \arrow{e}
\node{M(\ca)} \arrow{n,r}{~}
\end{diagram}
\]
We now define the trace morphism 
$\tau_Z: M(\per(\ca\on Z)) \ra \R\Ga_Z(X,M(\ca))$ as the
composition
\[
M(\per(\ca\on Z)) \ra \Ga(X, M(\roundper(\ca\on Z))) \ra
\R\Ga_Z(X,M(\ca)).
\]
We then have a commutative diagram
\[
\begin{diagram}
\node{M(\per(\ca\on Z))} \arrow{s} \arrow{e}
\node{M(\per \ca)} \arrow{s}\\
\node{\R\Ga_Z(X, M(\ca))} \arrow{e} 
\node{\R\Ga(X, M(\ca))}
\end{diagram}
\]
This yields a canonical lift of the classes 
constructed in section \ref{charclasses}
to the theories supported in $Z$. The trace
morphism $\tau_Z$ is invertible if $X$ and
$U=X\setminus Z$ are quasi-compact separated
schemes (by \ref{agreement} below).

\section{The main theorem, examples, proof}

This section is devoted to the main theorem
\ref{agreement}. Let $k$ be a field and $X$ a quasi-compact
separated scheme over $k$. The mixed complex associated with
$X$ is defined as $M(X)=\R\Ga(X,M(\co_X))$.
The main theorem states that the
trace map $\tau : M(\per X) \ra M(X)$ of
\ref{charclasses} is invertible in the
mixed derived category. 

In \ref{deffunc}, we define $M(\per X)$ and examine 
its functoriality with respect to morphisms
of schemes following \cite{Thomason90}. In \ref{agreement},
we state the theorem and, as a corollary, the case of
quasi-projective varieties. As an application,
we compute, in \ref{tiltingbundles},
the cyclic homology of smooth projective varieties
admitting a tilting bundle as described in the introduction.

The proof of the main theorem occupies subsections 
\ref{proofaffine} to \ref{proofofagreement}.
It proceeds by induction on the
number of open affines needed to cover $X$. The case of
affine $X$ is treated in section \ref{proofaffine}. 
The induction step uses a 
Mayer-Vietoris theorem (\ref{mayervietoris}) which is based
on the description of the fiber of the morphism of mixed
complexes induced by the localization at a quasi-compact
open subscheme. This description is achieved in
\ref{localperfect}. It is based on 
Thomason-Trobaugh's localization theorem, 
which we recall in section \ref{recallthomason}
in a suitable form, and on the localization theorem
for cyclic homology of localization pairs \cite[2.4 c)]{Keller97},
which we adapt to our needs in \ref{localizationcyclic}.

\subsection{Definition and functoriality.}
\label{deffunc} We adapt ideas of
Thomason-Trobaugh~\cite{Thomason90}: 
Let $X$ be a quasi-compact separated scheme over
a field $k$. We put 
$\per X=\per \co_X$ (cf. \ref{presheafperfect}). 
We claim that the assignment $X \mapsto M(\per X)$ is
a functor of $X$. Indeed, let $\flatper X$ be the pair
formed by the category of right bounded perfect complexes
with flat components and its subcategory of acyclic 
complexes. Then the inclusion
\[
\flatper X \ra \per X
\]
induces an equivalence in the associated triangulated
categories (by \cite[3.5]{Thomason90}) and hence an
isomorphism
\[
M(\flatper X) \ra M(\per X)
\]
by \cite[2.4 b)]{Keller97}. Now if $f: X \ra Y$ is a
morphism of schemes, then $f^*$ clearly induces a 
a functor $\flatper Y \ra \flatper X$ and hence
a morphism $M(\per Y) \ra M(\per X)$. Notice that
this morphism is compatible with the
map $M(\per X) \ra \R\Ga(X,M(\roundper X))$
of section \ref{charclasses}. 

Now suppose that $X$ admits an ample family of line
bundles. Then the inclusion
\[
\strper X \ra \per X
\]
induces an equivalence in the associated triangulated
categories \cite[3.8.3]{Thomason90} and hence an
isomorphism $M(\strper X) \ra M(\per X)$. Note that
$\strper X$ is simply the category of bounded complexes
over the category $\vecbundle X$ of algebraic vector bundles
on $X$ (together with its subcategory of acyclic complexes).
Hence we have the equality $M(\strper X)=M(\vecbundle X)$ where
$M(\vecbundle X)$ denotes the mixed complex associated with
the {\em exact category $\vecbundle X$} as defined in 
\cite{Keller97}.
In particular, if $X=\spec A$ is affine, we have canonical
isomorphisms
\[
M(A) \iso M(\proj A) \iso M(\vecbundle X) \iso M(\per X).
\]

\subsection{The main theorem.}
\label{agreement}
Let $X$ be a quasi-compact separated
scheme over a field $k$. The mixed complex associated with
$X$ is defined as $M(X)=\R\Ga(X,M(\co_X))$. Note that by
definition, we have
\[
\HCmixs(X)=\HCs M(X) \ko \quad \HNmixs(X) = \HNmixs M(X)\ko \ldots .
\]
\begin{theorem}  The trace morphism (\ref{charclasses})
\[
\tau : M(\per X) \ra M(X)
\]
is invertible. More generally, if $Z$ is a closed subset
of $X$ such that $U=X\setminus Z$ is quasi-compact, then
the trace morphism
\[
\tau_Z : M(\per (X \on Z)) \ra \R\Ga_Z(X, M(\co_X))
\]
is invertible.
\end{theorem}

\begin{corollary} Let $X$ be a quasi-compact separated
scheme over a field $k$. Then there is a canonical
isomorphism 
\[
\HCs(\per X)\iso \HCs(X).
\]
In particular, if $X$ admits an ample line bundle (e.g. if
$X$ is a quasi-projective variety), there is a canonical
isomorphism
\[
\HCders(\vecbundle X) \iso \HCs(X).
\]
\end{corollary}

The corollary was announced in \cite[1.10]{Keller97}, where
we wrote $\HCs(\vecbundle X)$ instead of $\HCders(\vecbundle X)$. It
is immediate from the theorem once we prove that for
quasi-compact separated schemes, there is an isomorphism
\[
\HCs(X) \iso \HCmixs(X).
\]
This will be done in \ref{proofofcorollary}.

The theorem will be proved in \ref{proofofagreement}.
The plan of the proof is described in the introduction
to this section.

\subsection{The example of varieties with tilting bundles.}
\label{tiltingbundles}
Suppose that $k$ is an algebraically
closed field and that $X$ is a smooth projective algebraic variety.
Suppose moreover that $X$ admits a {\em tilting bundle}, i.e. 
a vector bundle $T$ without higher selfextensions whose direct
summands generate the bounded derived category of the category
of coherent sheaves on $X$ as a triangulated category. 
Examples of varieties satisfying
these hypotheses are projective spaces,
Grassmannians, and smooth quadrics \cite{Beilinson78}, 
\cite{Kapranov83}, \cite{Kapranov86}, \cite{Kapranov88inv}.

\begin{proposition} The Chern character induces an isomorphism
\[
K_0(X)\ten_\Z \HCs(k) \ra \HCs(X).
\]
\end{proposition}
Here the left hand side is explicitly
known since the group $K_0(X)$ is free and admits a basis
consisting of the classes of the pairwise non-isomorphic
indecomposable direct summands of the tilting bundle.
For example, if $X$ is the Grassmannian of $k$-dimensional
subspaces of an $n$-dimensional space, the indecomposables
are indexed by all Young diagrams with at most $k$ rows
and at most $n-k$ columns.
Cyclic homology of projective spaces was first computed
by Beckmann \cite{Beckmann91} using a different method. 

The proposition shows that if $X$ is a smooth
projective variety such that $\H{n}(X,\co_X)\neq 0$
for some $n>0$, then $X$ cannot admit a tilting bundle.
Indeed, the group $\H{n}(X,\co_X)$ occurs as a direct
factor of $\HC{-n}(X)$ and therefore has to vanish
if the assumptions of the proposition are satisfied.

\begin{proof} Let $A$ be the endomorphism algebra of the
tilting bundle $T$ and $r$ the Jacobson radical of $A$. We 
assume without restriction of generality that $T$ is a 
direct sum of pairwise non-isomorphic indecomposable
bundles. Then $A/r$ is a product of copies of $k$
(since $k$ is algebraically closed). We will show
that the mixed complex $M(X)$ is canonically
isomorphic to $M(A/r)$. 
For this, consider the exact functor
\[
?\ten_A : \proj(A) \ra \vecbundle(X).
\]
It induces an equivalence in the bounded derived categories
\[
\cd^b(\proj(A)) \ra \cd^b(\vecbundle(X)).
\]
Indeed, we have a commutative square
\[
\begin{diagram}
\node{\cd^b(\proj(A))} \arrow{e,t}{?\ten_A T} \arrow{s}
\node{\cd^b(\vecbundle(X))} \arrow{s}
\\
\node{\cd^b(\mod(A))} \arrow{e,t}{\L(?\ten_A T)}
\node{\cd^b(\coh(X))} 
\end{diagram}
\]
where $\mod(A)$ denotes the abelian category of all
finitely generated right $A$-modules and $\coh(X)$
the abelian category of all coherent sheaves on $X$.
Since $T$ is a tilting bundle,
the bottom arrow is an equivalence.
Since $X$ is smooth projective, it follows that
$A$ is of finite global dimension. Hence the left
vertical arrow is an equivalence. Again because $X$ is
smooth projective, the right vertical arrow is an 
equivalence. Hence the top arrow is an equivalence.
So the functor
\[
?\ten_A T : \per(A) \ra \per(X)
\]
induces an equivalence in the associated triangulated
categories and hence an isomorphism
\[
M(\per(A)) \iso M(\per(X))
\]
by \cite[2.4 b)]{Keller97}. Of course, it also induces
an isomorphism $K_0(\proj(A)) \iso K_0(\vecbundle(X))$ and
the Chern character is compatible with these isomorphisms
by its description in \ref{charclasses}. So we are
reduced to proving that the Chern character induces
an isomorphism
\[
K_0(A) \ten_Z \HCs(k) \iso \HCs(A).
\]
For this, let $E\subset A$ 
be a semi-simple subalgebra such that $E$ identifies
with the quotient $A/r$. The algebra $E$ is a product
of copies of $k$ and of course, the inclusion $E\subset A$ 
induces an isomorphism in $K_0$. It also induces an
isomorphism in $\HCs$ by \cite[2.5]{Keller98}
since $A$ is finite-dimensional and of finite global
dimension. These isomorphisms
are clearly compatible with the Chern character and
we are reduced to the corresponding assertion for
$\HCs(E)$. This is clear since $E$ is a product
of copies of $k$.
\end{proof}

\subsection{Proof of the main theorem in the affine case.}
\label{proofaffine}
Suppose that $X=\spec A$. Then we know by section \ref{agreement}
that the canonical morphism
$M(A) \ra M(\per X)$ is invertible. Now Weibel-Geller
have shown in \cite[4.1]{WeibelGeller91}
that the canonical morphism
\[
M(A) \ra \R_{ce}\Ga(X, M(\co_X))
\]
is invertible where $M(\co_X)$ is viewed as a complex of sheaves
on $X$ and $\R_{ce}\Ga(X,?)$ denotes 
Cartan-Eilenberg hypercohomology (cf. section \ref{sheaveswith}).
Moreover, Weibel-Geller have shown in \cite[0.4]{WeibelGeller91} 
that the complex $M(\co_X)$ has quasi-coherent homology. By 
section \ref{sheaveswith}, it follows that the canonical morphism
\[
\R\Ga(X, M(\co_X)) \ra \R_{ce}\Ga(X, M(\co_X))
\]
is invertible. Using the commutative diagram
\[
\begin{diagram}
\node{M(A)} \arrow{s}\arrow{e} 
\node{\R_{ce}\Ga(X, M(\co_X))} 
\\
\node{M(\per X)} \arrow{e} 
\node{\R\Ga(X, M(\co_X))} \arrow{n}
\end{diagram}
\]
we conclude that $M(\per X) \ra\R\Ga(X, M(\co_X))$ is invertible
for affine $X$.

\subsection{Thomason-Trobaugh's localization theorem.}
\label{recallthomason}
Let $X$ be a quasi-compact quasi-separated scheme. We denote
by $\ct\per X$ the full subcategory of the (unbounded) derived
category of the category of $\co_X$-modules whose objects are the
perfect complexes. This category identifies with
the triangulated category associated with the
localization pair $\per X$ as defined in
\cite[2.4]{Keller97}. Recall that a triangle functor
$\cs\ra\ct$ is an {\em equivalence up to factors} if
it is an equivalence onto a full subcategory whose
closure under forming direct summands is all of $\ct$.
A sequence of triangulated categories
\[
0\ra \cR \ra \cs \ra \ct \ra 0
\]
is {\em exact up to factors} if the first functor
is an equivalence up to factors onto the kernel of
the second functor and the induced functor
$\cs/\cR\ra\ct$ is an equivalence up to factors.

\begin{theorem} {\em \cite{Thomason90}}
\begin{itemize}
\item[a)] Let $U\subset X$ be a quasi-compact open subscheme
and let $Z=X\setminus U$. Then the sequence
\[
0 \ra \ct\per(X \on Z) \ra \ct\per X \ra \ct\per U \ra 0
\]
is exact up to factors.
\item[b)] Suppose that $X=V\cup W$, where $V$ and $W$ are
quasi-compact open subschemes and put $Z=X\setminus W$. 
Then the lines of the diagram 
\[
\dgARROWLENGTH=1.5em
%\dgHORIZPAD=0.5em
\begin{diagram}
\node{\myspace 0} \arrow{e}
\node{\ct\per(X \on Z)} \arrow{s,l}{j^*} \arrow{e}
\node{\ct\per X} \arrow{s}\arrow{e}
\node{\ct\per W} \arrow{s}\arrow{e}
\node{0\myspace } \\
\node{\myspace 0} \arrow{e}
\node{\ct\per(V \on Z)} \arrow{e}
\node{\ct\per V} \arrow{e}
\node{\ct\per(V\cap W)} \arrow{e}
\node{0\myspace}
\end{diagram}
\]
are exact up to factors and the functor $j^*$ is an equivalence
up to factors.
\end{itemize}
\end{theorem}

The theorem was proved in section 5 of \cite{Thomason90}.
Note that the first assertion of part b) follows from a).
The second assertion of b) is a special case of the
main assertion in \cite[5.2]{Thomason90} (take 
$U=V$, $Z=X\setminus W$ in [loc.cit.]). 
A new proof of the theorem is due to  A. Neeman \cite{Neeman92},
\cite{Neeman96}.

\subsection{Localization in cyclic homology of DG categories.}
\label{localizationcyclic}
In this section, we adapt the localization theorem
\cite[4.9]{Keller97} to our needs. Let 
\[
0 \ra \ca \arr{F}  \cb \arr{G} \cc \ra 0
\]
be a sequence of small flat exact DG categories such that
$F$ is fully faithful, $GF=0$, and the induced sequence of
stable categories
\[
0 \ra \ul{\ca}\ra \ul{\cb} \ra \ul{\cc} \ra 0
\]
is exact up to factors (\ref{recallthomason}). 
\begin{theorem} The morphism
\[
\cone(M(\ca) \xrightarrow{M(F)} M(\cb)) \ra M(\cc)
\]
induced by $M(G)$ is a quasi-isomorphism.
\end{theorem}

\begin{proof} The proof consists in extracting the relevant
information from \cite{Keller97} : Indeed, since $F$ is fully
faithful, we may consider $\ca\arr{F}\cb$ as a localization
pair and since $GF=0$, the square
\[
\begin{diagram}
\node{\ca} \arrow{s}\arrow{e,t}{F} 
\node{\cb} \arrow{s,r}{G} 
\\
\node{0} \arrow{e}
\node{\cc}
\end{diagram}
\]
as a morphism of localization pairs, i.e. a morphism of the
category $\cl^b_{str}$ of \cite[4.3]{Keller97}. By applying the
completion functor $?^+$ of [loc. cit.] we obtain a morphism
\begin{equation}
\begin{diagram}
\node{(\ca \arr{F} \cb)^+} \arrow{s} \\
\node{(0 \ra \cc)^+}
\end{diagram}
\label{themorphism}
\end{equation}
of the category $\cl$. Applying the functor $Cm$ to this
morphism yields the morphism
\[
\begin{diagram}
\node{(M(\ca) \ra M(\cb))} \arrow{s,r}{(0,M(G))}
\\
\node{(0 \ra M(\cc))}
\end{diagram}
\]
of $\dmormix$ by the remarks following proposition
4.3 of \cite{Keller97}. On the other hand, applying
the functor $I_\lambda$ of \cite[4.8]{Keller97} to the
morphism (\ref{themorphism}) yields the identity
of $\cc^+$ in $\cm$ and applying $M$ (denoted by $C$
in \cite{Keller97}) yields the identity of $M(\cc)$ in
$\dmix$. By the naturality of the isomorphism of
functors in \cite[4.9 a)]{Keller97}, call it $\psi$, 
we obtain a commutative square in $\dmix$
\[
\begin{diagram}
\node{\cone(M(\ca) \ra M(\cb))} \arrow{e,t}{\psi} \arrow{s,l}{(0,M(G))}
\node{M(\cc)}  \arrow{s,r}{\id}
\\
\node{\cone(0 \ra M(\cc))} \arrow{e,t}{\psi}
\node{M(\cc)}
\end{diagram}
\]
So the left vertical arrow of the square is 
invertible in $\dmix$, which is what we had to prove.
\end{proof}

\subsection{Perfect complexes with support and local cohomology.}
\label{localperfect}
Let $X$ be a quasi-compact quasi-separated scheme, $U\subset X$
a quasi-compact open subscheme, and $Z=X\setminus U$.
Let $j: U \ra X$ be the inclusion.

\begin{proposition} The sequence
\[
M(\roundper(X\on Z)) \ra M(\roundper X) \ra j^* M(\roundper U)
\]
embeds into a triangle of $\dmix(X)$. This triangle is
canonically isomorphic to the $Z$-local cohomology triangle
associated with $M(\roundper X)$. In particular, there is
a canonical isomorphism
\[
M(\roundper (X\on Z)) \iso \R\Ga_Z(X,M(\roundper X)).
\]
Moreover, the canonical morphisms fit into a morphism of triangles
\[
\dgARROWLENGTH=1.5em
\begin{diagram}
\node{M(\per(X \on Z))} \arrow{s}\arrow{e}
\node{M(\per X)} \arrow{s}\arrow{e}
\node{M(\per U)} \arrow{s}\arrow{e}
\node{M(\per(X \on Z))[1]} \arrow{s}
\\
\node{\Ga_Z M(\roundper X))} \arrow{e}
\node{\Ga M(\roundper X))} \arrow{e}
\node{\Ga M(\roundper U))} \arrow{e}
\node{\Ga_Z M(\roundper X))[1]} 
\end{diagram}
\]
in the mixed derived category, where $\Ga$ and $\Ga_Z$
are short for $\R\Ga(X,?)$ and $\R\Ga_Z(X,?)$.
\end{proposition}

\begin{proof} Let $V\subset X$ be open. Consider the
sequence
\begin{equation}
M(\per(V\on Z)) \ra M(\per V) \ra M(\per (V\cap U)).
\label{crudesequence}
\end{equation}
If we let $V$ vary, it becomes a sequence of presheaves on $X$.
We will show that there is a sequence of mixed complexes
of presheaves
\begin{equation}
A \arr{f} B \arr{g} C
\label{finesequence}
\end{equation}
such that 
\begin{itemize}
\item we have $gf=0$ in the category of mixed complexes of presheaves
\item in the derived category of mixed complexes of presheaves, the
sequence \ref{finesequence} becomes isomorphic to the
sequence \ref{crudesequence}.
\item for each quasi-compact open subscheme $V\subset X$, the canonical
morphism from the cone over the morphism $A(V) \ra B(V)$
to $C(V)$ induced by $g$ is a quasi-isomorphism.
\end{itemize}
This implies that firstly, the sequence of sheaves associated with
the sequence \ref{finesequence} embeds canonically into a triangle
\[
\tilde{A} \ra \tilde{B} \ra \tilde{C} \ra \tilde{A}[1] \ko
\]
where the tilde denotes sheafification and the connecting
morphism is constructed as the composition
\[
\tilde{C} \liso \cone(\tilde{A} \ra \tilde{B}) \ra \tilde{A}[1]\ko
\]
and secondly we have a morphism of triangles
\[
\begin{diagram}
\node{A(V)} \arrow{s} \arrow{e}
\node{B(V)} \arrow{s} \arrow{e}
\node{C(V)} \arrow{s} \arrow{e}
\node{A(V)[1]} \arrow{s} 
\\
\node{\R\Ga(V, \tilde{A})} \arrow{e}
\node{\R\Ga(V, \tilde{B})} \arrow{e}
\node{\R\Ga(V, \tilde{C})} \arrow{e}
\node{\R\Ga(V, \tilde{A}[1])}
\end{diagram}
\]
for each quasi-compact open subscheme $V\subset X$ (to prove this
last assertion, we use that $\R\Ga(V,?)$ lifts to a derived
functor defined on the category of all sequences 
\[
A' \arr{f'} B' \arr{g'} C'
\]
with $g' f'=0$).  

To construct the sequence \ref{finesequence}, we have to (pre-) sheafify
a part of the proof of \cite[2.4]{Keller97}. For this, let  $\iper X$ denote
the category of all fibrant (\ref{terminology})
perfect complexes. Then the
inclusion $\iper X \ra \per X$ induces an equivalence in the
associated triangulated categories and thus we have an 
isomorphism $M(\iper X) \iso M(\per X)$ in $\dmix$.  Note that
this even holds if $X$ is an arbitrary ringed space. In particular,
it holds for each open subscheme $V\subset X$ instead of $X$.
Hence the presheaf $V \mapsto M(\per V)$ is isomorphic
in the derived category of presheaves to $V \mapsto M(\iper V)$.
Similarly for the other terms of the sequence, so that we
are reduced to proving the assertion for the sequence
of presheaves whose value at $V$ is
\[
M(\iper (V \on Z)) \ra M(\iper V) \ra M(\iper (U\cap V)).
\]
For this, let $\ci(V)$ be the exact dg category \cite[2.1]{Keller97}
of fibrant (\ref{terminology}) complexes on $V$ and let 
$\tilde{\ci}(V)$ be the category whose objects are
the exact sequences 
\[
0 \ra K \arr{i} L \arr{p} M \ra 0
\]
of $\ci(V)$ such that $i$ has split monomorphic components,
$K$ is acyclic off $Z$ and $i_x$ is a quasi-isomorphism for each
$x\in Z$. Then $\tilde{\ci}(V)$ is equivalent to a full exact dg
subcategory of the category of filtered objects of $\ci(V)$
(cf. example 2.2 d) of \cite{Keller97}). Let $\tilde{\ci}(V\on Z)$
be the full subcategory of $\tilde{\ci}(X)$ whose objects are the
sequences 
\[
0 \ra K \iso L \ra 0 \ra 0
\]
and $\tilde{\ci}(U\cap V)$ the full subcategory whose objects
are the sequences
\[
0 \ra 0 \ra M \iso L \ra 0.
\]
Let $G : \tilde{\ci}(V) \ra \tilde{\ci}(V\cap U)$ be the functor
\[
(0 \ra K \ra L \ra M \ra 0) \mapsto (0 \ra 0 \ra M \arr{\id} M \ra 0)
\]
and $F:\tilde{\ci}(V\on Z) \ra \tilde{\ci}(V)$ the inclusion. Then
the sequence
\begin{equation}
0 \ra \tilde{\ci}(V\on Z) \arr{F} \tilde{\ci}(V) \arr{G} 
\tilde{\ci}(V\cap U) \ra 0 
\end{equation}
is an exact sequence of the category $\cm_{str}$ of
\cite[4.4]{Keller97} and in particular we have $GF=0$. 
We take the subsequence of perfect objects :
Let $\tiper (V \on Z)$ be the full subcategory of $\tilde{\ci}(V\on Z)$
whose objects are the $K \iso L \ra 0$ with $K\in\iper(V\on Z)$,
let $\tiper(V)$ be the full subcategory of the $K \ra L \ra M$
with $M \in \per V$, and let $\tiper(V\cap U)$ be the full
subcategory of the $0 \ra L \ra M$ with $M|_U \in \per(V\cap U)$.
Consider the diagram
\[
\dgARROWLENGTH=2em
\begin{diagram}
\node{\tiper(V\on Z)} \arrow{e,t}{F} \arrow{s}
\node{\tiper(V)} \arrow{e,t}{G} \arrow{s}
\node{\tiper(V\cap U)} \arrow{s}
\\
\node{\iper(V\on Z)} \arrow{e}
\node{\iper(V)} \arrow{e}
\node{\iper(V\cap U),}
\end{diagram}
%\begin{CD}
%\tiper(V\on Z) @>F>> \tiper(V) @>G>> \tiper(V\cap U) \\
%    @VVV              @VVV             @VVV          \\
%\iper(V\on Z)  @>>>  \iper(V)  @>>>   \iper(V\cap U)
%\end{CD}
\]
where the three vertical functors are given by
\begin{eqnarray*}
K \iso L \ra 0 & \mapsto & K \\
K \ra L \ra M  & \mapsto & M \\
0 \ra L \ra M  & \mapsto & M|_U .
\end{eqnarray*}
Its left hand square is commutative up to isomorphism and its 
right hand square is commutative up to the homotopy 
\cite[3.3]{Keller97}
\[
L|_U \xrightarrow{p|_U} M|_U .
\]
The vertical arrows clearly induce equivalences in the 
associated triangulated categories. By applying the functor
$M$ to the diagram and letting $V$ vary we obtain a
commutative diagram in the derived category of presheaves of mixed
complexes on $X$. The vertical arrows become invertible
and the top row becomes
\[
M(\tiper(V\on Z)) \ra M(\tiper(V)) \ra M(\tiper(V\cap U))
\]
where $V$ runs through the open subsets of $X$. This
is the sequence of presheaves $A \ra B \ra C$ announced at 
the beginning of the proof. Using 
theorem \ref{recallthomason} a) and theorem \ref{localizationcyclic}
one sees that it has the required properties.
\end{proof}

\subsection{Mayer-Vietoris sequences.}
\label{mayervietoris} Let $X$ be a quasi-compact quasi-separated
scheme and $V,W\subset X$ quasi-compact open subschemes
such that $X=V\cup W$.

\begin{proposition} There is a canonical morphism of triangles
in the mixed derived category
{
\overfullrule=0pt
$$
\dgARROWLENGTH=1.5em
\dgHORIZPAD=0.5em
\begin{diagram}
\node{M(\per X)} \arrow{s}\arrow{e}
\node{M(\per V)\oplus M(\per W)} \arrow{s} \arrow{e}
\node{M(\per (V\cap W) )} \arrow{s}\arrow{e}
\node{\mybigspace}\\
\node{\Ga M(\roundper X)} \arrow{e}
\node{\Ga M(\roundper V)\oplus \Ga M(\roundper W))} \arrow{e}
\node{\Ga M(\roundper (V\cap W))} \arrow{e}
\node{\mybigspace}
\end{diagram}
$$
}
\overfullrule=5pt
where $\Ga$ is short for $\R\Ga(X,?)$.
\end{proposition}

\begin{proof}
Put $Z=X\setminus W$. 
The first line of the diagram is deduced from theorem \ref{recallthomason} b)
using \cite[2.7]{Keller97}. Clearly the two squares appearing in
the diagram are commutative. We have to show that the
square involving
the arrows of degree $1$ 
\[
\begin{diagram}
\node{M(\per (V\cap W) )} \arrow{s}\arrow{e}
\node{M(\per X)[1]} \arrow{s} \\
\node{\Ga M(\roundper (V\cap W))} \arrow{e}
\node{\Ga M(\roundper X))[1]}
\end{diagram}
\]
is commutative as well. By [loc.cit.], 
the connecting morphism is the composition
\[
\dgARROWLENGTH=1.5em
\begin{diagram}
\node{}
\node{M(\per(X \on Z)[1]} \arrow{s}\arrow{e}
\node{M(\per X)[1]}
\\
\node{M(\per(V\cap W))} \arrow{e}
\node{M(\per(V\on Z))[1]} 
\node{}
\end{diagram}
\]
Here the vertical morphism is invertible by theorem \ref{recallthomason} b)
and \cite[2.4 b)]{Keller97}. The second line of the diagram is
the Mayer-Vietoris triangle for hypercohomology. So the connecting
morphism of the second line is obtained as the composition
\[
\Ga M(\roundper(V\cap W))\ra
\Ga_{Z}M(\roundper V)[1] \liso
\Ga_{Z}M(\roundper X)[1] \ra
\Ga M(\roundper X)[1] \ko
\]
where $\Ga$ and $\Ga_Z$ are short for $\R\Ga(X,?)$ and $\R\Ga_Z(X,?)$.
Now it follows from proposition \ref{localperfect} that the rightmost square
of the diagram of the assertion is commutative as well.
\end{proof}

\subsection{Proof of theorem 5.1.}
\label{proofofagreement}
Let $V_1, \ldots, V_n$ be open affines covering $X$. 
If $n=1$, theorem \ref{agreement}
holds by section \ref{proofaffine}. If $n>1$, we cover
$X$ by $V=V_1$ and $W=\bigcup_{i=2\ldots n} V_i$.
The intersection $V \cap W$ is then covered by 
the $n-1$ sets $V \cap V_i$, $2\leq i \leq n$. These are affine,
since $X$ is separated. So theorem \ref{agreement}
holds for $V$, $W$, and $V\cap W$ by the induction
hypothesis. Thus it holds for $X=V\cup W$ by
proposition \ref{mayervietoris}. 
The assertion for $\tau_Z$ now follows by
proposition \ref{localperfect}.

\subsection{Proof of corollary 5.1.}
\label{proofofcorollary} In \cite{Weibel96} (cf. also \cite{Weibel97}),
C.~Weibel defined $HC_*(X)$ as the homology of the complex of
$k$-modules
\[
\R\Gamma_{ce}(X, CC(\co_X))
\]
where $\R\Gamma_{ce}$ denotes Cartan-Eilenberg hypercohomology
(cf. section \ref{sheaveswith})
and $CC(\co_X)$ is the sheafification of the 
classical bicomplex. Now Weibel-Geller have
shown in \cite{WeibelGeller91} that the 
Hochschild complex $C(\co_X)$ has
quasi-coherent homology. Thus each column of $CC(\co_X)$
has quasi-coherent homology and hence (the sum total
complex of) $CC(\co_X)$ has itself quasi-coherent homology.
Hence by theorem \ref{sheaveswith}, the above complex is 
isomorphic to
\[
\R\Gamma (X, CC(\co_X)).
\]
Now, as in the case of an algebra (cf. \cite[2.5.13]{Loday92}),
$CC(\co_X)$ may also be viewed as the (sum total complex of the)
bicomplex $\cb C(M(\co_X))$ associated with the mixed complex of
sheaves $M(\co_X)$ (cf. \cite[Section 2]{Weibel97}). 
What remains to be proved then
is that the canonical map
\[
\cb C(\R\Gamma(X, M(\co_X)) \ra \R\Gamma(X, \cb C(M(\co_X))
\]
is invertible in the derived category of $k$-vector spaces.
Now indeed, more generally, we claim that we have
\[
\cb C(\R\Gamma(X, M)) \iso \R\Gamma(X, \cb C(M))
\]
for any mixed complex of sheaves $M$ with quasi-coherent
homology. As the reader will easily check, this is immediate 
once we know that the functor $\R\Gamma(X, ?)$ commutes with
countable direct sums when restricted to the category
of complexes with quasi-coherent homology. This follows
from Corollary 3.9.3.2 in \cite{Lipmanxy}. It may
also be proved by the argument of \cite[1.4]{Neeman96}. For
completeness, we include a proof : 
Let $K_i$, $i\in I$, be a family of complexes with quasi-coherent
homology. It is enough to prove that $H^0(X,?)$ takes
$K=\bigoplus K_i$ to the sum of the $H^0(X, K_i)$. 
Now $\Gamma(X,?)$ is of finite cohomological dimension
on the category of quasi-coherent modules.
Indeed, for an affine $X$, this follows
from Serre's theorem \cite[III, 1.3.1]{GrothendieckDieudonne61},
and for arbitrary $X$ it is proved by induction on
the size of an affine cover of $X$ (here we use that 
$X$ is quasi-compact and separated). It therefore follows from 
by theorem \ref{sheaveswith} b), lemma \ref{acyclicimage},
and Serre's theorem \cite[III, 1.3.1]{GrothendieckDieudonne61}.
that we have an isomorphism 
$H^0(X, K_i) \iso H^0(X, \tau^{\,\geq n} K_i)$
and similarly for $K$ for some fixed $n<0$
(cf. the proof of theorem \ref{sheaveswith} for the
definition of the truncation functor $\tau^{\,\geq n}$).
So we may assume
that the $K_i$ and $K$ are uniformly bounded below. But
then, we may compute the $H^0(X,K_i)$ using resolutions 
$K_i \ra F_i$ by uniformly bounded below complexes of 
flasque sheaves.
The sum of the $F_i$ is again bounded below with
flasque components and is clearly quasi-isomorphic to $K$. 
Now $\Ga(X,?)$ commutes with infinite sums since $X$
is quasi-compact, so the claim follows.

%%%%%%%%%%%%%%%%%%%%%%%%%%%%%%%%%%%%%%%%%%%%%%%%%%%%%%%%%%%%%%%%%%%%
\appendix
%%%%%%%%%%%%%%%%%%%%%%%%%%%%%%%%%%%%%%%%%%%%%%%%%%%%%%%%%%%%%%%%%%%%

\section{On Cartan-Eilenberg resolutions}

We prove that Cartan-Eilenberg hypercohomology coincides with
derived functor hypercohomology on all (unbounded) complexes of sheaves
with quasi-coherent homology on a quasi-compact separated scheme.
More precisely, we prove that in this situation, Cartan-Eilenberg
resolutions are actually $K$-injective resolutions in the sense of
\cite{Spaltenstein88}. 

\subsection{Terminology.}
\label{terminology}

Let $\ca$ be a Grothendieck category.
Spaltenstein \cite{Spaltenstein88} defined a complex 
$I$ over $\ca$ to be {\em $K$-injective} if, in the homotopy category,
there are no non zero morphisms from an acyclic complex
to $I$. This is the case iff each morphism $M \ra I$
in the derived category is represented by a unique
homotopy class of morphisms of complexes.

In \cite[A.2]{Weibel97}, C.~Weibel proposed the use of the term {\em
fibrant} for $K$-injective. Indeed, one can show that a complex is
$K$-injective iff it is homotopy equivalent to a complex which is
fibrant for the `global' closed model structure on the category of 
complexes in which cofibrations are the componentwise monomorphisms.
This structure is an additive analogue of the global closed
model structure on the category of simplicial sheaves on
a Grothendieck site. The existence of the global structure
in the case of simplicial sheaves was proved by
Joyal \cite{Joyal84} (cf. \cite[2.7]{Jardine87}).
We have not been able to find a published proof 
of the fact that the category of complexes over a 
Grothendieck category admits the global structure 
(an unpublished proof is due to F.~Morel). 
However, the key step may be found in \cite[Prop. 1]{Frankexy}.

Whereas in the homotopy category, the notions
of `fibrant for the global structure' and `$K$-injective' 
become essentially equivalent, there is a slight difference 
at the level of complexes: fibrant objects for the global 
structure are exactly the $K$-injective 
complexes with injective components.

We will adopt the terminology proposed
by Weibel: We {\em call a complex fibrant}
iff it is $K$-injective in the sense of Spaltenstein.
This will not lead to ambiguities since we will 
not use the global closed model structure.

\subsection{Sheaves with quasi-coherent cohomology.}
\label{sheaveswith}
Let $X$ be a scheme and $K$ a complex of $\co_X$-modules
(unbounded to the right and to the left). Let $I$
be a Cartan-Eilenberg resolution of $K$, i.e.
\begin{itemize}
\item[a)] $I$ is a $\Z\times\Z$-graded $\co_X$-module endowed
with differentials $d_I$ of bidegree $(1,0)$ and $d_{II}$ of
bidegree $(0,1)$ such that $(d_I+d_{II})^2=0$,
\item[b)] $I^{pq}$ vanishes for $q<0$ and
\item[c)] $I$ is endowed with an augmentation $\eps : K \ra I$, 
i.e. a morphism of differential $\Z\times\Z$-graded $\co_X$-modules,
where $K$ is viewed as concentrated on the $p$-axis, such that for each
$p$, the induced morphisms $K^p \ra I^{p,\bt}$ and
$H^p K \ra H^p_I I$ are injective resolutions. 
\end{itemize}
It follows that for each $p$, the induced morphisms
$B^p K \ra B^p_I K$ and $Z^p K \ra Z^p_I I$ are injective resolutions
and that the rows of $I$ are products of complexes of the form
\[
\ldots 0 \ra M \ra 0 \ldots \quad\mbox{or} 
\quad \ldots 0 \ra M \arr{\id} M \ra 0 \ldots \ko
\]
where $M$ is injective.

Let $J=\hTot I$ denote the product total complex of $I$ and
$\eta : K \ra J$ the morphism of complexes induced by $\eps$.
The morphism $\eta$ is called a {\em total Cartan-Eilenberg resolution}
of $K$. The {\em Cartan-Eilenberg hypercohomology} of $K$ is the
cohomology of the complex
\[
\R\Ga_{ce}(X,K) = \Ga(X, J).
\]
The morphism $\eta$ is usually {\em not} a quasi-isomorphism.

\begin{theorem} 
\begin{itemize}
\item[a)] The complex $J$ is fibrant (\ref{terminology}).
\item[b)] If $K$ has quasi-coherent homology, the morphism
$\eta : K \ra J$ is a quasi-isomorphism. Hence,
Cartan-Eilenberg hypercohomology of $K$ coincides with
derived functor hypercohomology of $K$ in the sense of 
Spaltenstein \cite{Spaltenstein88}.
\end{itemize}
\end{theorem}

Part a) holds more generally whenever $K$ is a complex of objects over
an abelian category having enough injectives and admitting all
countable products. This was proved by C.~Weibel in
\cite[A.3]{Weibel96}.  For completeness, we include a proof 
of a) below. Part b) was proved by C. Weibel in [loc. cit.] for the case of
complete abelian categories with enough injectives and 
{\em exact products}, for example module categories.
The case we consider here is implicit in \cite[3.13]{Spaltenstein88}.
Nevertheless, we thought it useful to include the explicit statement
and a complete proof.

In preparation of the proof, let us recall the notion of a homotopy limit
(cf. \cite{BoekstedtNeeman93} for example) : If $\ct$
is a triangulated category admitting all countable products
and 
\[
 \ldots \ra X_{p+1} \arr{f_p} X_{p} \ra \ldots \ra X_0 \ko p\in\N \ko
\]
is a sequence in $\ct$, its {\em homotopy limit} 
$\holim X_p$ is defined by the Milnor triangle \cite{Milnor62}
\[
\holim X_p \ra \prod X_p \arr{\Phi} \prod X_q \ra (\holim X_p)[1] \ko
\]
where the morphism $\Phi$ has the components
\[
\prod_p X_p \xrightarrow{\can} X_{q+1}\oplus X_q
\xrightarrow{[-f_{q}, \id]} X_q.
\]
Note that the homotopy limit is unique only up to non
unique isomorphism. We will encounter the following
situation : Consider a sequence of complexes
\[
\ldots \ra K_p \arr{f_p} K_{p-1} \ra \ldots \ra K_0
\]
over an additive category admitting all countable products
such that the $f_p$ are componentwise split epi
(or, more generally, for each $n$ and $p$, the morphism 
$X^n_{p+k} \ra X^n_p$ is split epi for some $k\gg 0$). Then we have
a componentwise split short exact sequence of complexes
\[
0 \ra \lii K_p \ra \prod_p K_p \arr{\Phi} \prod_q K_q \ra 0
\]
and hence the inverse limit $\lii K_p$ is then isomorphic
to $\holim K_p$ in the homotopy category.

\begin{proof}[Proof of the theorem]
a) Note that the bicomplex $I$ is the inverse limit of its
quotient complexes $I^{\bt,q]}$ obtained by killing all rows of
index greater than $q$. 
Let $J_q$ be the product total complex of $I^{\bt,q]}$.
Then the sequence of the $J_q$ has inverse limit $\hTot I$
and its structure maps are split epi in each component.
Hence $I$ is isomorphic to the homotopy limit of the sequence
of the $J_q$.
Since the class of fibrant complexes
is stable under extensions and products, it is stable
under homotopy limits. Therefore it is 
enough to show that the $J_q$ are fibrant.
Clearly the $J_q$ are iterated extensions of
rows of $I$ (suitably shifted). So it is enough to
show that the rows of $I$ are fibrant.
But each row of $I$ is homotopy equivalent to a complex
with vanishing differential and injective components.
Such a complex is the product of its components placed
in their respective degrees and is thus fibrant.

b) For $p\in\Z$, define $\tau^{\,\geq p}K$ to be the quotient
complex of $K$ given by
\[
 \ldots \ra 0 \ra K^p/B^p K \ra K^{p+1} \ra K^{p+2} \ra \ldots
\]
and $\tau^{\, <p}K$ to be the subcomplex of $K$ given by
\[
\ldots \ra K^{p-2} \ra K^{p-1} \ra B^p K \ra 0 \ra \ldots .
\]
Define $\tau^{\,\geq p} J$ and $\tau^{\, <p}J$ by applying the
respective functor to each row of $J$. Then the 
morphism $\tau^{\,\geq p} K \ra \tau^{\,\geq p}J$ is a 
Cartan-Eilenberg resolution for each $p\in\Z$. Since
$\tau^{\,\geq p}K$ is left bounded, it follows that the
induced morphism $\tau^{\,\geq p}K \ra \hTot \tau^{\,\geq p} J$
is a quasi-isomorphism for each $p\in\Z$. Now fix
$n\in\Z$ and consider the diagram
\[
\begin{diagram}
\node{H^n K} \arrow{e} \arrow{s}
\node{H^n \tau^{\,\geq p} K} \arrow{s} \\
\node{H^n \hTot J} \arrow{e} \node{H^n \hTot \tau^{\,\geq p} J .}
\end{diagram}
\]
For $p<n$, the top morphism is invertible. It now
suffices to show that for $p\ll 0$, the bottom morphism
is invertible. Equivalently, it is enough to show that
$H^n \hTot \tau^{\, <p} J$ vanishes for $p\ll 0$. For
this let $x\in X$. We have to show that $(H^n \hTot \tau^{\, <p} J)_x$
vanishes. Since taking the stalk is an exact functor,
this reduces to showing that the complex
$(\hTot \tau^{\, <p} J)_x$ is acyclic in degree $n$. For
this, it is enough to show that $(\hTot \tau^{\, <p} J)(U)$
is acyclic in degree $n$ for each affine neighbourhood of $x$.
Now $\tau^{\, <p} J$ is a Cartan-Eilenberg resolution
of $\tau^{\, <p} K$. Therefore, if we apply proposition
\ref{acyclicimage} below to the functor
$F=\Gamma(U,?)$, we see that $(\hTot \tau^{\, <p} J)(U)$
is acyclic in all degrees $n\geq p$. Indeed, we have
$(\R^{i} F)(H^p K)=0$ for all $p$ and all $i>0$ by
Serre's theorem \cite[III, 1.3.1]{GrothendieckDieudonne61}, 
since $H^p K$ is quasi-coherent.

\end{proof}

\subsection{Unbounded complexes with uniformly bounded cohomology.}
\label{acyclicimage}
Let $\ca$ be an abelian category with enough injectives
which admits all countable products and let 
$F: \ca \ra \ca b$ be an additive functor commuting with
all countable products. 

Let $K$ be a complex over $\ca$
and let $K \ra J$ a Cartan-Eilenberg resolution.
Suppose that $K^p=0$ for all $p>0$ and 
that there is an integer $n$ with
\[
(\R^i F)(H^p K)=0
\]
for all $i\geq n$ and all $p\in\Z$. 

\begin{lemma} 
We have $H^p F\hTot J=0$ for all $p\geq n$.
\end{lemma}

Note that this assertion is clear if
$K$ is (homologically) left bounded. The point is that
it remains true without this hypothesis.

\begin{proof} Define $\tau^{\,\geq p} K$ and $\tau^{\,\geq p} J$
as in the proof of proposition \ref{sheaveswith}.
The canonical morphisms $\tau^{\,\geq p} J \ra \tau^{\,\geq p+1} J$
are split epi in each bidegree and $J$ identifies
with the inverse limit of the $\tau^{\,\geq p} J$. Hence
we have $\hTot J = \lii \hTot \tau^{\,\geq p}J$ and
the morphisms 
\[
\hTot \tau^{\,\geq p} J \ra \hTot \tau^{\,\geq p+1} J
\]
are componentwise split epi. Since $F$ commutes with
countable products, we therefore have
$F (\hTot J) = \lii F \hTot \tau^{\,\geq p} J$. 
By lemma \ref{mitlef} below, it is therefore enough
to show that the groups $H^i F(\hTot L_p)$ vanish
for all $i\geq n$ and all $p$ where $L_p$ is the
kernel of the canonical morphism 
$\tau^{\,\geq p} J \ra \tau^{\,\geq p+1} J$. Now $L_p$ is in fact
a Cartan-Eilenberg resolution of the kernel of
the morphism $\tau^{\,\geq p} K \ra \tau^{\,\geq p+1} K$, which
is isomorphic to the complex
\[
\ldots 0 \ra K^{p-1}/ B^{p-1} K \ra Z^p K \ra 0 \ra \ldots
\]
This complex is quasi-isomorphic to $H^p K$ placed in degree $p$.
So $\hTot L_p$ is homotopy equivalent to an injective
resolution of $H^p K$ shifted by $p$ degrees. Hence
\[
H^i F \hTot L_p = H^i \R F(H^p K [-p]) = (\R^{i-p} F)(H^p K).
\]
By assumption, this vanishes for $i-p\geq n$. 
\end{proof}

\subsection{A Mittag-Leffler lemma.}
\label{mitlef}
Let $n$ be an integer and let
\[
\ldots \ra K_{p+1} \longarr{\pi_{p+1}} K_p \ra \ldots \ra 
K_0 \arr{\pi_0} K_{-1}=0 \ko p\in\N \ko
\]
be an inverse system of complexes of abelian groups such
that the $\pi_p$ are surjective in each component and
$H^i K'_p=0$ for all $i\geq n$ and all $p$, where
$K'_p$ is the kernel of $\pi_{p}$.

\begin{lemma}
We have $H^i \lii K_p=0$ for all $i\geq n$.
\end{lemma}

\begin{proof}
By induction, we find that $H^i K_p=0$ for all $i\geq n$.
Now we have exact sequences
\[
0 \ra Z^i K_p \ra K^i_p \ra Z^{i+1}K_p \ra 0 \ko
\]
for all $i\geq n-1$. Since $B^i K_p \iso Z^i K_p$, the
maps $Z^i K_{p+1} \ra Z^i K_p$ are surjective for $i \geq n$.
The fact that $H^n K'_{p+1}=0$ implies that the
maps $Z^{n-1} K_{p+1} \ra Z^{n-1} K_p$ are surjective as well.
By the Mittag-Leffler lemma 
\cite[$0_{III}$, 13.1]{GrothendieckDieudonne61},
the sequence
\[
0 \ra \lii Z^i K_p \ra \lii K^i_p \ra \lii Z^{i+1} K_p \ra 0
\]
is still exact for $i\geq n-1$. 
Since $\lii Z^i K_p \liso Z^i \lii K_p$,
this means that $H^i \lii K_p=0$ for $i\geq n$.
\end{proof}

\section{A comparison of derived categories}

\subsection{Boekstedt-Neeman's theorem.}
\label{quasicoherentsheaves}
Let $X$ be a quasi-compact separated sche\-me, $\cd \qcoh X$
the derived category of the category $\qcoh X$ of quasi-coherent
sheaves on $X$, $\cd X$ the derived category of all sheaves
of $\co_X$-modules on $X$, and $\cd_{qc}X$ its full subcategory
whose objects are the complexes with quasi-coherent homology.

As an application of theorem \ref{sheaveswith}, we give
a partially new proof of the following result of Boekstedt-Neeman. 
We refer to \cite[Prop. 1.3]{AlonsoJeremiasLipman97} for yet
another proof.

\begin{theorem} {\em \cite[5.5]{BoekstedtNeeman93}}
The canonical functor $\cd\qcoh X \ra \cd_{qc}X$
is an equivalence of categories.
\end{theorem}

The proof proceeds by induction on
the size of an affine cover of $X$. The crucial step
is the case where $X$ is affine. Our proof for this case
is new. For completeness, we have included the full
induction argument.

\begin{proof}  In a first step, suppose that
$X$ is affine : $X = \spec A$. We identify $\qcoh X$ with
$\Mod A$ and then have to show that the sheafification functor
$F : \cd \Mod A \ra \cd X$ induces an equivalence
$\cd\Mod A \ra \cd_{qc} X $. Clearly, the image of $A$
(viewed as a complex of $A$-modules concentrated in degree $0$)
is $\co_X$. By the lemma below, it
suffices therefore to show that
\begin{itemize}
\item[a)] We have $A \iso \HOm{\cd X}{\co_X}{\co_X}$
and $\HOm{\cd X}{\co_X}{\co_X[n]}=0$ for each $n \neq 0$,
\item[b)] The object $\co_X$ is compact in $\cd_{qc} X$
i.e. the associated functor 
\[
\HOm{\cd_{qc} X}{\co_X}{?}
\]
commutes with infinite direct sums.
\item[c)] An object $K \in \cd_{qc} X$ vanishes if
$\HOm{\cd X}{\co_X}{K[n]}$ vanishes for all $n\in\Z$.
\end{itemize}
The three assertions a), b), and c) all follow easily
from the fact that we have an isomorphism
\[
\HOm{\cd_{qc} X}{\co_X}{?} \iso \Gamma(X, H^0(?)) \ko
\]
which we will now prove : Indeed, let $K\in \cd_{qc} X$.
By definition, we have 
\[
\HOm{\cd X}{\co_X}{K} = H^0 \R\Gamma(X,K).
\]
Now we have morphisms
\[
H^0 \R\Gamma(X,K) \larr{\alpha} H^0\R\Gamma(X,\tau_{\leq 0} K)
\arr{\beta} H^0\R\Gamma(X, H^0 K) = \Gamma(X, H^0 K).
\]
The morphism $\alpha$ is invertible because $\R\Gamma(X,?)$
is a right derived functor. The morphism $\beta$ is invertible
by theorem \ref{sheaveswith} b), lemma \ref{acyclicimage},
and Serre's theorem \cite[III, 1.3.1]{GrothendieckDieudonne61}.

Now suppose that $X$ is the union of
$n$ open affine sets $U_1, \ldots, U_n$. By induction on
$n$ and the affine case, we may assume that the claim
is proved for $U=U_1$ and $V=\bigcup_{i=2\ldots n} U_i$.
Let $j_1 : U \ra X$ and $j_2 : V \ra X$ be the inclusions.
Let $Y=X\setminus U$ and let $i: Y \ra X$ be the inclusion.
For any object $K\in \cd_{qc} X$, we have a triangle
\[
\R\Gamma_Y K \ra K \ra j_{1*}j_1^* K \ra \R\Gamma_Y K[1].
\]
Here the second morphism is the adjunction morphism and
$\R\Gamma_Y K$ is defined (up to unique isomorphism) by
the triangle. The object 
$j_1^*K$ is a complex of sheaves on $U$ and
$H^n j_1^* K= j_1^* H^n K$ is quasi-coherent. So 
$j_1^* K$ is in the faithful image of $\cd\qcoh U$. Because
$X$ is separated, $j_{1*}$ preserves quasi-coherence
(cf. \cite[3.9.2]{Lipmanxy}). So the triangle lies
in $\cd_{qc} X$. The subset $Y\subset X$
is a closed subset of $V$ and $i=j_1 \, i_2$, where
$i_2$ is the inclusion of $Y$ into $V$. This implies that
$\R\Gamma_Y K = j_{2*} (\R\Gamma_{Y\subset V} K)$. 
The above triangle thus shows that
$\cd_{qc} X$ is generated by the 
$j_{1*}K'$ and the $j_{2*}K''$, where $K'$ belongs
to $\cd\qcoh U$ and $K''$ to $\cd\qcoh V$. It remains
to be checked that morphisms between 
$j_{1*}K'$ and $j_{2*}K''$ in $\cd\Mod\co_X$ are in
bijection with those in $\cd\qcoh X$. Indeed, we have
\[
\HOm{\cd X}{j_{1*}K'}{j_{2*}K''} =
\HOm{\cd V}{j_2^* j_{1*}K'}{K''}.
\]
By the induction hypothesis, the latter group identifies
with
\[
\HOm{\cd\qcoh V}{j_2^* j_{1*} K'}{K''} =
\HOm{\cd\qcoh X}{j_{1*}K'}{j_{2*}K''}.
\]
The same argument applies to morphisms from 
$j_{2*}K''$ to $j_{1*} K'$. This ends the proof.
\end{proof}

\subsection{Derived categories of modules.}
Let $A$ be a ring and $\ct$ a triangulated category admitting
all (infinite) direct sums. Suppose that $F: \cd\Mod A \ra \ct$ 
is a triangle functor commuting with all direct sums.
For the convenience of the reader, we include a proof of the
following more and more well-known 

\begin{lemma} The functor $F$ is an equivalence if and only if
\begin{itemize}
\item[a)] We have $A \iso \HOm{\ct}{FA}{FA}$ and
$\HOm{\ct}{FA}{FA[n]}=0$ for all $n\neq 0$.
\item[b)] The object $FA$ is compact in $\ct$, i.e.
$\HOm{\ct}{FA}{?}$ commutes with infinite direct sums.
\item[c)] An object $X$ of $\ct$ vanishes iff $\HOm{\ct}{FA}{X[n]}=0$
for all $n\in\Z$.
\end{itemize}
\end{lemma}

\begin{proof} Let $\cs\subset\ct$ be the smallest triangulated
subcategory of $\ct$ containing $FA$ and stable under
forming infinite direct sums. Then, since $FA$ is compact,
the inclusion $\cs\ra\ct$ admits a right
adjoint $R$ by Brown's representability theorem 
\cite{Brown62} (cf. also \cite[5.2]{Keller94}, 
\cite{Neeman96}, \cite{Frankexy}). 
Now if $X\in\ct$ and $RX \ra X \ra X' \ra RX[1]$
is a triangle over the adjunction morphism, then
$\HOm{\ct}{FA}{X'[n]}$ vanishes for all $n\in\Z$
by the long exact sequence associated with the
triangle. So $X'$ vanishes by assumption c) and
$\cs$ coincides with $\ct$. So $FA$ is a compact
generator for $\ct$. Now the claim follows from
\cite[4.2]{Keller94}.
\end{proof}

\Addresses


\begin{thebibliography}{TT}
\normalsize

\bibitem{AvramovFoxbyHalperinxy} L.~Avramov, H.-B.~Foxby, 
S.~Halperin, {\em Differential graded homological algebra},
preprint in preparation.

\bibitem{Beckmann91} C.~Beckmann, {\em Zyklische Homologie
von Schemata und relative algebraische $K$-Theorie},
Inaugural-Dissertation, K\"oln, 1991.

\bibitem{Beilinson78} A.~A.~Beilinson, {\em Coherent sheaves on
{\nobf P}$^n$ and problems of linear algebra}, Funkts. Anal. Prilozh.
{\nobf 12} (1978), 68-69. English translation: Funct. Anal.
Appl. {\nobf 12} (1979), 214-216.

\bibitem{BoekstedtNeeman93} M.~Boekstedt, A.~Neeman,
{\em Homotopy limits in triangulated categories},
Comp. Math. {\nobf 86} (1993), 209--234.

\bibitem{BresslerNestTsygan97}
P.~Bressler, R.~Nest, B.~Tsygan, {\em Riemann-Roch Theorems
via Deformation Quantization}, Preprint, 
alg-geom/9705014, June 1997.

\bibitem{Brown62}  E. H. Brown, {\em Cohomology theories},
Ann. of Math. {\nobf 75} (1962), 467-484.

\bibitem{Deligne73} P.~Deligne, {\em Cohomologie \`a supports propres},
Expos\'e XVII, SGA 4, Springer LNM {\nobf 305} (1973), 252-480.

\bibitem{Frankexy} J.~Franke, {\em On the Brown representability
theorem for triangulated categories}, Preprint, 1997.

\bibitem{GrothendieckDieudonne61}
A.~Grothendieck, J.~Dieudonn\'e, {\em El\'ements de g\'eom\'etrie
alg\'ebrique}, Publ Math. I.H.E.S. No. 11, (1961), 5--167.

\bibitem{Jardine87} J.~F.~Jardine, {\em Simplicial presheaves},
J. Pure Applied Alg. {\nobf 47} (1987), 35-87.

\bibitem{Joyal84} A.~Joyal, {\em Letter to Grothendieck}, 1984.

\bibitem{Kapranov83} M.~M.~Kapranov, {\em The derived categories of
coherent sheaves on Grassmannians}, Funkts. Anal. Prilozh.
{\nobf 17} (1983), 78-79. English translation: Funct. Anal.
Appl. {\nobf 17} (1983), 145-146.

\bibitem{Kapranov86} M.~M.~Kapranov, {\em The derived category of 
coherent sheaves on a quadric}, Funkts. Anal. Prilozh. {\nobf
20} (1986), 67. English Translation: Funct. Anal. Appl.
{\nobf 20} (1986), 141--142.

\bibitem{Kapranov88inv} M.~M.~Kapranov, {\em On the derived categories
of coherent sheaves on some homogeneous spaces}, Invent.
math. {\nobf 92} (1988), 479--508.

% \bibitem{Kapranov88izv} M.~M.~Kapranov, {\em On the derived category
% and $K$-functor of coherent sheaves on intersections of quadrics},
% Izv. Akad. Nauk SSSR Ser. Matem. {\nobf 52} (1988), 186--199.
% English Translation: Math. USSR Izvestiya {\nobf 32} (1989),
% 191--204.

\bibitem{Kassel87} C.~Kassel, {\em Cyclic homology, comodules 
and mixed complexes}, J.~Alg. {\nobf 107} (1987), 195--216.

\bibitem{Keller94} B.~Keller, 
{\em Deriving DG categories}, Ann.~scient.~Ec.~Norm.~Sup., 4$^e$ s\'erie
{\nobf 27} (1994), 63--102.

\bibitem{Keller98} B.~Keller,
{\em Invariance and Localization for cyclic homology of
DG algebras}, J. Pure Appl. Alg. {\nobf 123} (1998), 223--273.

\bibitem{Keller97} B.~Keller, {On the cyclic homology of exact
categories}, Preprint Nr. 185 in www.math.uiuc.edu/K-theory/,
to appear in J. Pure Appl. Alg.


\bibitem{Lipmanxy} J.~Lipman, {\em Notes on derived categories
and derived functors}, Preprint 
(available at http://www.math.purdue.edu/$\tilde{\;\;\;}$lipman/).

\bibitem{Loday86} J.-L.~Loday, {\em Cyclic homology: a survey},
Banach Center Publications {\nobf 18} (1986), 285--307.

\bibitem{Loday92} J.-L.~Loday, {\em Cyclic homology},
Grundlehren 301, Springer-Verlag Berlin Heidelberg, 1992.

\bibitem{McCarthy94} R.~McCarthy, 
{\em The cyclic homology of an exact category}, 
J. Pure and Appl. Alg. {\nobf 93} (1994), 251--296.

\bibitem{Milnor62} J.~Milnor, {\em On Axiomatic Homology Theory},
Pacific J.~Math. {\nobf 12} (1962), 337--341.

\bibitem{Mitchell72} B.~Mitchell, {\em Rings with several objects},
Adv. in Math. {\nobf 8} (1972), 1--161.

\bibitem{Neeman92} A.~Neeman, {\em The Connection between the
$K$-theory localization theorem of Thomason, Trobaugh and
Yao and the smashing subcategories of Bousfield and
Ravenel}, Ann.~Scient.~\'Ec.~Norm.~Sup. {\nobf 25} (1992),
547--566.

\bibitem{Neeman96} A.~Neeman, {\em The Grothendieck duality 
theorem via Bousfield's techniques and Brown representability}, 
J. Am. Math. Soc. {\nobf 9} (1996), 205--236.

\bibitem{Quillen73} D.~Quillen, {\em Higher Algebraic $K$-theory I},
Springer LNM {\nobf 341}, 1973, 85--147.

\bibitem{SchapiraSchneiders94} P.~Schapira, J.-P.~Schneiders,
{\em Index theorem for elliptic pairs, Elliptic pairs II},
Ast\'erisque vol. 224, 1994.

\bibitem{Spaltenstein88} N. Spaltenstein, 
{\em Resolutions of unbounded complexes}, 
Compositio Mathematica {\nobf 65} (1988), 121--154.

\bibitem{AlonsoJeremiasLipman97} L.~Alonso~Tarr\'{\i}o, 
A.~Jerem\'{\i}as L\'opez, J.~Lipman, 
{\em Local Homology and Cohomology on Schemes},
Ann. scient. \'Ec. Norm. Sup., 4e s\'erie, t. 30, 1997, 1--39.

\bibitem{Thomason90} 
R.~W.~Thomason, T.~Trobaugh, {\em Higher Algebraic 
$K$-Theory of Schemes and of Derived Categories}, in {\em The
Grothendieck Festschrift}, III, Birkh\"auser, 
Progress in Mathematics {\nobf 87} (1990), p. 247--436.

\bibitem{Weibel96} C.~A.~Weibel, 
{\em Cyclic Homology For Schemes}, 
Proc. Amer. Math. Soc. {\nobf 124} (1996), 1655--1662.

\bibitem{Weibel97} C.~A.~Weibel,
{\em The Hodge Filtration and Cyclic Homology},
K-Theory {\nobf 12} (1997), 145--164.

\bibitem{WeibelGeller91} C.~A.~Weibel, S.~C.~Geller,
{\em \'Etale descent for Hochschild and cyclic homology},
Comment. Math. Helv. {\nobf 66} (1991), 368--388.

\end{thebibliography}
\end{document}